\documentclass[11pt, leqno]{article}

\usepackage{hyperref}
\usepackage{array} 
\usepackage{amsfonts}
\usepackage{amsmath}
\usepackage{amssymb}
\usepackage{graphicx}
\usepackage{color}
\usepackage{esint}    

\usepackage{enumitem} 

\newtheorem{thm}{Theorem}[section]
\newtheorem{cor}[thm]{Corollary}
\newtheorem{lem}[thm]{Lemma}
\newtheorem{prop}[thm]{Proposition}
\newtheorem{defn}[thm]{Definition}
\newtheorem{rem}[thm]{Remark}
\numberwithin{equation}{section}

\newcommand{\dx}{\,{\rm d}x}
\newcommand{\dy}{\,{\rm d}y}

\newcommand{\dt}{\,{\rm d}t}
\newcommand{\rd}{{\rm d}}

\def\LL{\mathrm{L}} 

\newcommand{\A}{\mathcal{L}}
\newcommand{\AM}{\mathcal{L}^{\frac{1}{2}}}
\newcommand{\AI}{\mathcal{L}^{-1}}
\newcommand{\AIM}{\mathcal{L}^{-\frac{1}{2}}}

\newcommand{\n}{F}

\newcommand{\K}{{\mathbb K}}

\newcommand{\RR}{\mathbb{R}}
\newcommand{\R}{\mathbb{R}}

\def\dist{\mathrm{dist}} 

\def\qed{\,\unskip\kern 6pt \penalty 500
\raise -2pt\hbox{\vrule \vbox to8pt{\hrule width 6pt
\vfill\hrule}\vrule}\par}

\definecolor{darkblue}{rgb}{0.05, .05, .65}
\definecolor{darkgreen}{rgb}{0.1, .65, .1}
\definecolor{darkred}{rgb}{0.8,0,0}

\begin{document}
\title{\textbf{Infinite speed of propagation and regularity of solutions \\ to the fractional porous medium equation \\in general domains\\[7mm]}}

\author{Matteo Bonforte$^{\,a}$, Alessio Figalli$^{\,b}$, Xavier Ros-Oton$^{\,c}$}
\date{} 

\maketitle

\thispagestyle{empty}

\begin{abstract}
We study the positivity and regularity of solutions to the fractional porous medium equations $u_t+(-\Delta)^su^m=0$ in $(0,\infty)\times\Omega$, for $m>1$ and $s\in (0,1)$ and with Dirichlet boundary data $u=0$ in $(0,\infty)\times(\R^N\setminus\Omega)$, and nonnegative initial condition $u(0,\cdot)=u_0\geq0$.

Our first result is a quantitative lower bound for solutions which holds for all positive times $t>0$.
As a consequence, we find a global Harnack principle stating that for any $t>0$ solutions are comparable to $d^{s/m}$, where $d$ is the distance to $\partial\Omega$.
This is in sharp contrast with the local case $s=1$, where the equation has finite speed of propagation.

After this, we study the regularity of solutions.
We prove that solutions are classical in the interior ($C^\infty$ in $x$ and $C^{1,\alpha}$ in $t$) and establish a sharp $C^{s/m}_x$ regularity estimate up to the boundary.

Our methods are quite general, and can be applied to wider classes of nonlocal parabolic equations of the form $u_t+\mathcal L F(u)=0$ in $\Omega$, both in bounded or unbounded domains.
\end{abstract}

\vfill

\vspace{.5mm}
\begin{itemize}
\item[(a)] Departamento de Matem\'{a}ticas, Universidad Aut\'{o}noma de Madrid,\\
Campus de Cantoblanco, 28049 Madrid, Spain.\\
E-mail:\texttt{~matteo.bonforte@uam.es}.\\
Web-page:\texttt{~http://www.uam.es/matteo.bonforte}
\item[(b)] Department of Mathematics, The University of Texas at Austin,\\
2515 Speedway Stop C1200, Austin, TX 78712-1082, USA.\\
E-mail:\texttt{~figalli@math.utexas.edu}.\\
Web-page:\texttt{~http://www.ma.utexas.edu/users/figalli/}
\item[(c)] Department of Mathematics, The University of Texas at Austin,\\
2515 Speedway Stop C1200, Austin, TX 78712-1082, USA.\\
E-mail:\texttt{~ros.oton@math.utexas.edu}.\\
Web-page:\texttt{~http://www.ma.utexas.edu/users/ros.oton/}
\end{itemize}

\small
\tableofcontents
\normalsize

\section{Introduction}

Our goal is to investigate the regularity properties of nonnegative solutions to nonlinear nonlocal diffusion equations of degenerate type that can be written in the form
\begin{equation}\label{NNL.EQ}
u_t=-\A \n (u).
\end{equation}
Here $\A$ is a linear operator possibly of fractional order, and $\n$ is a monotone nondecreasing function, satisfying some conditions that allow for degeneracies, like $\n(0)=0$ and $F'(0)=0$. To see (and understand) the effect of degeneracy, we impose accordingly zero Dirichlet boundary conditions. The prototype equation is given by the so-called \textsl{Fractional Porous Medium Equation}
\begin{equation}\label{FPME.EQ}
u_t+(-\Delta)^{s}u^m=0\quad \textrm{in}\quad (0,\infty)\times\Omega,
\end{equation}
with Dirichlet conditions $u\equiv0$ in $(0,\infty)\times(\RR^N\setminus\Omega)$.
Here $m>1$ and $s\in(0,1)$, $\Omega\subset \RR^N$ is a bounded $C^{1,1}$ domain, and $\A=(-\Delta)^{s}$ is the Fractional Laplacian, namely
\[(-\Delta)^{s}f(x)=c_{N,s}{\rm PV}\int_{\RR^N} \frac{f(x)-f(x+y)}{|y|^{N+2s}}\dy,\]
where $c_{N,s}>0$ is an explicit constant whose exact value is not relevant to our purposes.
We have decided to concentrate here in this particular model, but most of our techniques can be extended to more general equations of the form
\eqref{NNL.EQ} with the operator $\A$ given  by
\begin{equation}
\label{eq:general L}\A f(x):={\rm PV}\int_{\RR^N}\big(f(x)-f(x+y)\big)\,\K(x,y)\dy
\end{equation}
under appropriate conditions on the kernel $\K$ and on the inverse $\AI$ of the operator $\A$.
Also, the boundedness of $\Omega$ is not essential as most of our estimates are of local nature.
We postpone the discussion to these more general situations to Section \ref{sec.generalizations}, where further details are given, together with a discussion about the possible extension of our results also to more general nonlinearities $\n$.

Nonlinear nonlocal diffusion models of this type have received a lot of attention in the last years,
especially because of their applications to anomalous diffusions in physics and biology; see \cite{BGT2000,LMT2003,SV2013}.
These equations appear also as hydrodynamic limits of interacting particle systems with long-range dynamics, cf. \cite{Jara0,Jara1,JKOlla,Jara2}, and also in boundary heat control problems \cite{AC2010,DL1972}. We refer the interested reader to \cite{AC2010,BV2012,BV-PPR1,DPQRV2,SV2013,VazAbel,VazSurvey} for further details about possible applications.

In this paper we are going to prove positivity, global H\"older regularity, and interior higher regularity results, for a general class of nonnegative weak solutions, called weak dual solutions; this class has been introduced in \cite{BV-PPR1} and contains various other classes of solutions such as mild (semigroup), weak, weak energy, or $H^{-s}$ solutions, as discussed in \cite{BV-PPR1,BV-PPR2-1,BSV2013}. In Appendix I, we recall the definition of weak dual solutions and their basic properties, together with the relation with other notions of solutions.

One of the main results of this paper is to show that nonnegative weak dual solutions are indeed classical in the interior of the domain: we prove that, inside $\Omega$, they are $C^\infty_{\rm loc}$ in space and $C^{1,\alpha}_{\rm loc}$ in time.
Furthermore, we study the boundary regularity of solutions and establish a $C^{\frac{s}{m},\frac{1}{2m}}_{x,t}$ H\"older estimate up to the boundary, which is optimal in the $x$ variable.

Let us mention that, while the case $s=1$ (corresponding to the classical Porous Medium Equation) has been extensively studied in the last 30 years by many authors \cite{Ar69,Ar70b,ArCaVa85,CaFr79b, CaFr80,CaVaWo87,CaWo90,DaHaLe01,DaHa98,DiB86,DiB88,DGVacta,GiPe81}, see also the books \cite{DaskaBook, DiBook,DGVbook,VazBook}, not many results are currently known when $s \in (0,1)$.
Indeed, to the best of our knowledge, the only two results available in this setting are the papers of Athanasopoulos and Caffarelli \cite{AC2010},
where the authors prove the global $C^\alpha$ regularity  of solutions in space-time, and the one of Vazquez, dePablo, Quiros and Rodriguez \cite{VDPQR}, where it is shown that weak solutions of the Cauchy problem posed on the whole space are classical.
We point out that some arguments in the latter paper crucially exploits the fact that the problem is set on the whole space,
as they use a fractional version of the Aleksandrov reflection principle / moving planes method.
Also, the paper of Athanasopoulos and Caffarelli relies on the so-called ``extension property'' for the fractional Laplacian,
and does not generalize to general operators as in \eqref{eq:general L}.
Hence, new ideas and techniques have to be introduced in our situation.

A crucial issue addressed in this paper is the positivity of solutions, which opens the road to our regularity results.
Our quantitative and precise lower bounds hold for all times whenever $s\in (0,1)$, and  reveal a peculiar property of the nonlocal evolution ($0<s<1$) versus the local one ($s=1$): indeed this shows that solutions have infinite speed of propagation when $0<s<1$,
in sharp contrast with the local case $s=1$. More comments on this important issue will be given below.

\noindent\textit{Positivity and Harnack type bounds. }Under suitable assumptions on the inverse of the linear operator~$\A$, there is a quite complete theory of existence and uniqueness of a quite general class of weak solutions, as well as a priori estimates, cf. \cite{BV-PPR1, BV-PPR2-1, BV-PPR2-2}. In those papers the Global Harnack Principle (GHP) has been proved to hold after a waiting time $t_*$ for all $s \in (0,1]$:
\[
h_0 \frac{\Phi_1(x)^{\frac{1}{m}}}{t^{\frac{1}{m-1}}} \le u(t,x)\le h_1\frac{\Phi_1(x)^{\frac{1}{m}}}{t^{\frac{1}{m-1}}}\,,\qquad\mbox{for all $t\ge t_*$ and all $x\in \overline{\Omega}$\,,}
\]
with some explicit positive constants $h_0,h_1$. The GHP turns out to imply a localized version of Harnack inequalities that hold at the same time and for the same balls, and are more similar to the classical Harnack inequalities:
\[
\sup_{x\in B_R(x_0)} u(t,x)\le \mathcal{H}\inf_{x\in B_R(x_0)} u(t,x)\,,
\qquad\mbox{for all $t\ge t_*$\,.}
\]
It is however remarkable the fact that we can take the same time in both members of the Harnack inequality.
Actually, even backward (in time) Harnack inequalities hold true; see Appendix I for more details.  We recall that Fabes, Garofalo, and Salsa \cite{FGS86} showed that backward Harnack inequalities hold for the case $s=m=1$; indeed, such backward inequalities can be true only for the Dirichlet problem on bounded domains and fail for solutions to the Cauchy problem on the whole space, as counterexamples show.

We notice that, while the upper bound in the GHP holds for any $t>0$, the lower bound (hence the local version) was only known to hold after a suitable waiting time $t_*$, that has been quantified in \cite{BV-PPR1, BV-PPR2-1, BV-PPR2-2} as follows:
\[
t_*(u_0)=\frac{c_*}{\|u_0\|_{\LL^1_{\Phi_1}(\Omega)}^{m-1}}\,,
\]
with explicit $c_*>0$ depending only on $N,s,m,\Omega$.

\noindent\textit{Finite versus infinite speed of propagation. }When $s=1$\,, i.e. when dealing with the classical (local) PME, this waiting time cannot be avoided, in view of the finite speed of propagation, see \cite{Ar-Pe, VazBook}.
On the other hand, when $s<1$, { it was conjectured in  \cite{BV-PPR1, DPQRV2} that the fractional PME should have infinite speed of propagation. }Here we answer positively to this question: in Theorem \ref{Thm1} below, we prove a quantitative lower bounds for all positive times, which clearly imply infinite speed of propagation. In the case of the Cauchy problem on $\RR^N$, the first qualitative positivity statements for all times can be found in \cite{DPQRV2}\,, and in a more quantitative way using lower barriers in Theorem 1.4  of \cite{SV2013}. The question of finite vs infinite speed of propagation on $\RR^N$ has been also investigated in \cite{CV2011,CV2010,SdTV2,SdTV3} for a different model of fractional PME introduced by Caffarelli and Vazquez \cite{CV2011,CV2010}; an equivalence among selfsimilar solutions  between the two main models of nonlocal porous medium equations has been established in \cite{SdTV1}, and this gives an hint about why we shall expect infinite speed of propagation also for the PME model under investigation.

\noindent\textit{Regularity. } As mentioned above, the quantitative positivity bounds are crucial for us to attack the issue of regularity of solutions.
Indeed, using these bounds we establish first the global $C^\alpha$ regularity of solutions in space and time, and then we adapt the arguments in \cite{VDPQR}
to obtain the $C^\infty_x$ regularity in the interior of the domain.
Furthermore, we also establish the optimal $C_x^{s/m}$ regularity \emph{up to the boundary} for solutions $u$ to \eqref{FPME.EQ}, and prove the H\"older regularity in $t$ up to the boundary.
Finally, this boundary regularity estimate in $t$ combined with the regularity in space, yields the interior $C^{1,\alpha}_t$ regularity.

\noindent\textbf{Plan of the paper. }In the following two subsections we state our main results on the positivity and regularity of solutions. Section 2 contains the proof of the positivity results, namely Theorem \ref{Thm1} and its consequences Theorems \ref{Thm2} and \ref{Thm3}; to justify our argument,
 proofs are performed first for solutions of a suitable approximate problem,
 and then we pass to the limit in the estimates. The existence and basic properties of these approximate solutions are collected in Appendix II.
Section 3 contains the proofs of the regularity results, Theorems \ref{Thm4}, \ref{Thm5}, and \ref{Thm6}.
Finally, we discuss in Section 4 how the techniques used in this paper can be extended to solutions of more general nonlocal diffusion equations of degenerate type, and which kind of results extend to those settings.
Appendix I contains a discussion on the precise functional analytic setting that we use, together with the definition and basic properties of weak dual solutions.

After the writing of this paper was completed, we learned that,
at the very same time and independently of us, de Pablo, Quir\'os and Rodr\'iguez proved
in  \cite{DPQR} the H\"older regularity for solutions to the Cauchy problem in $\R^N$ for equation of the form \eqref{NNL.EQ} with general kernels and nonlinearities, see Remark \ref{rem.PQRV} in Section \ref{sec.generalizations}.

\subsection{Main results on positivity and Harnack estimates.}

\noindent Consider the following homogeneous Dirichlet problem
\begin{equation}\label{Dirichlet.Problem}
\left\{
\begin{array}{lll}
\partial_t u+(-\Delta)^s u^m=0 & \qquad\mbox{for any $(t,x)\in (0,\infty)\times \Omega$}\\
u(t,x)=0 &\qquad\mbox{for any $(t,x)\in (0,\infty)\times (\RR^N\setminus\Omega)$}\\
u(0,x)=u_0(x) &\qquad\mbox{for any $x\in\Omega$}\,.
\end{array}
\right.
\end{equation}
Our first main result concerns quantitative lower bounds.
Here and in the following, we denote by $\Phi_1$ the first eigenfunction for the operator $\A$
(see Appendix I), and we use the notation
$$
\|u_0\|_{\LL^1_{\Phi_1}(\Omega)}:=\int_{\Omega}u_0\,\Phi_1\dx.
$$
We recall here that, from now on, we will always consider $m>1$, $s\in (0,1)$ and $N>2s$, unless explicitly stated. It is worth spending a few words about the condition $N>2s$. Our approach is based on estimates valid for weak dual solutions, a class of weak solutions corresponding to a dual problem, involving the inverse of the fractional Laplacian $\AI=(-\Delta)^{-s}$, which has a kernel given by the Green function, see Appendix I for further details. Notice that the Green function of the fractional Laplacian has a singularity $|x|^{-(N-2s)}$ when $N>2s$\,, a logarithmic singularity when $N=2s$ and stops to be singular when $N<2s$\,, cf. \cite{Chen-Song,Kul}.
The latter two cases appear only in dimension $N=1$ and deserve a different treatment, even if most of our techniques could be successfully adapted.

\begin{thm}[Global quantitative positivity]\label{Thm1}
Let $m>1$, $0<s<1$, and $N>2s$. Let $\Omega$ be a bounded domain of class $C^{1,1}$,
and let $u$ be a solution to the Dirichlet problem \eqref{Dirichlet.Problem} corresponding to a nonnegative initial datum $u_0\in \LL^1_{\Phi_1}(\Omega)$\,.
Then the following bound holds true:
\begin{equation}\label{Thm1.lower.bdd}
u(t,x)\ge \kappa\|u_0\|^m_{\LL^1_{\Phi_1}(\Omega)}\,t\,\Phi_1(x)^{\frac{1}{m}}\,,\qquad\mbox{for all $0\le t\le t_*$ and all $x\in \overline{\Omega}$\,,}
\end{equation}
where $t_*=t_*(u_0)$ has the form
\begin{equation}\label{t*.Thm1}
t_*(u_0)=\frac{c_*}{\|u_0\|_{\LL^1_{\Phi_1}(\Omega)}^{m-1}}\,,
\end{equation}
and $c_*,\kappa$ are positive constants depending only on $N,s,m,\Omega$. As a consequence, solutions to the Dirichlet problem \eqref{Dirichlet.Problem} corresponding to nonnegative and nontrivial initial data have \textsl{infinite speed of propagation.}
\end{thm}

\noindent\textbf{Remark. }We stress on the fact that the constants $c_*,\kappa>0$ have an explicit form (see \eqref{kappa} for $\kappa$ and  \cite{BV-PPR1, BV-PPR2-1, BV-PPR2-2} for $c_*$).

Combining the result above with the Global Harnack Principle of \cite{BV-PPR1, BV-PPR2-1, BV-PPR2-2}, we can immediately obtain the following two theorems.

\begin{thm}[Global Harnack Principle for all times]\label{Thm2}
Let $m>1$, $0<s<1$, and $N>2s$.
Let $\Omega$ be a bounded domain of class $C^{1,1}$, let $u$ be a solution to the Dirichlet problem \eqref{Dirichlet.Problem} corresponding to a nonnegative initial datum $u_0\in \LL^1_{\Phi_1}(\Omega)$, and let $t_*$ be as in \eqref{t*.Thm1}.
Then
\begin{equation}\label{Thm2.GHP}
\underline{\kappa}\left(1\wedge\frac{t}{t_*}\right)^{\frac{m}{m-1}} \frac{d(x)^{\frac{s}{m}}}{t^{\frac{1}{m-1}}} \le u(t,x)\le \overline{\kappa}\,\frac{d(x)^{\frac{s}{m}}}{t^{\frac{1}{m-1}}}\,,\qquad\mbox{for all $t>0$ and all $x\in \overline{\Omega}$\,,}
\end{equation}
where $d(x)={\rm dist}(x,\partial\Omega)$, and $\overline{\kappa},\underline{\kappa}>0$ depend only on $N,s,m,\Omega$.
\end{thm}

The constants $\overline{\kappa},\underline{\kappa}>0$ have an explicit form given by \eqref{kappa} for $\underline{\kappa}$, and in \cite{BV-PPR1, BV-PPR2-2} for $\overline{\kappa}$.

\begin{thm}[Local Harnack inequalities for all times]\label{Thm3}
Let $m>1$, $0<s<1$, and $N>2s$. Let $\Omega$ be a bounded domain of class $C^{1,1}$,
 let $u$ be a solution to the Dirichlet problem \eqref{Dirichlet.Problem} corresponding to a nonnegative initial datum $u_0\in \LL^1_{\Phi_1}(\Omega)$,
 and let $t_*$ be as in \eqref{t*.Thm1}.
Then, for all balls $B_R(x_0)\subset\subset \Omega$,
\begin{equation}\label{Thm3.local.Harnack}
\sup_{x\in B_R(x_0)} u(t,x)\le \frac{\mathcal{H}}{\left(1\wedge\frac{t}{t_*}\right)^{\frac{m}{m-1}}}\inf_{x\in B_R(x_0)} u(t,x)\,,
\qquad\mbox{for all $t>0$\,,}
\end{equation}
where $\mathcal{H}>0$ depends only on $N,s,m,\Omega,\dist(B_R(x_0),\partial\Omega)$.
\end{thm}

\subsection{Main results on higher and boundary regularity.}

After proving quantitative positivity and Harnack type estimates for arbitrary positive times, we study the regularity of solutions.
We prove first the interior $C^\alpha$ regularity in space and time, which combined with the results of \cite{VDPQR} lead to the interior $C^\infty$ regularity in  space.
Moreover, we also establish the sharp regularity up to the boundary in the $x$-variable, namely we show that $u(t,\cdot)\in C^{s/m}(\overline\Omega)$ for all $t>0$.
This is stated in the following:

\begin{thm}[H\"older regularity up to the boundary]\label{Thm4}
Let $m>1$, $0<s<1$, and $N>2s$.
Let $\Omega$ be a bounded domain of class $C^{1,1}$,
and let $u$ be a solution to the Dirichlet problem \eqref{Dirichlet.Problem} corresponding to a nonnegative initial datum $u_0\in \LL^1_{\Phi_1}(\Omega)$.
Then, for each $0<t_0<T$ we have
\begin{equation}\label{Thm4.bdry-reg}
\|u\|_{C^{\frac{s}{m},\,\frac{1}{2m}}_{x,t}\left(\overline\Omega\times [t_0,T]\right)}\leq C,
\end{equation}
where $C$ depends only on $N,s, m, \Omega,t_0$, and $\|u_0\|_{\LL^1_{\Phi_1}(\Omega)}$.
\end{thm}

Notice that the $C^{s/m}_x$ regularity up to the boundary is optimal.
Indeed, by Theorem \ref{Thm2} we have that $u\geq c(u_0,t)d^{s/m}$, with $c(u_0,t)>0$ for all $t>0$, and therefore $u(t,\cdot)\notin C^{\frac{s}{m}+\epsilon}_x(\overline\Omega)$ for any $\epsilon>0$.

The interior $C^\infty$ regularity of the solution in the $x$-variable is given by the following result.

\begin{thm}[Higher interior regularity in space]\label{Thm5}
Let $m>1$, $0<s<1$, and $N>2s$.
Let $\Omega$ be a bounded domain of class $C^{1,1}$, and let $u$ be a solution to the Dirichlet problem \eqref{Dirichlet.Problem} corresponding to a nonnegative initial datum $u_0\in \LL^1_{\Phi_1}(\Omega)$\,.
Then $u\in C^\infty_x((0,\infty)\times\Omega)$.

More precisely, let $k\geq1$ be any positive integer, and $d(x)={\rm dist}(x,\partial\Omega)$.
Then, for any $t\geq t_0>0$ we have
\begin{equation}\label{Thm5.C^infty-reg}
\bigl|D^k_x u(t,x)\bigr|\leq C\,[d(x)]^{\frac{s}{m}-k},
\end{equation}
where $C$ depends only on $N,s, m,k,\Omega,t_0$, and $\|u_0\|_{\LL^1_{\Phi_1}(\Omega)}$.
\end{thm}

Finally, we study the higher interior regularity in time, and we get the following result.

\begin{thm}[$C^{1,\alpha}$ interior regularity in time]\label{Thm6}
Let $m>1$, $0<s<1$, and $N> 2s$.
Let $\Omega$ be a bounded domain of class $C^{1,1}$, and let $u$ be a solution to the Dirichlet problem \eqref{Dirichlet.Problem} corresponding to a nonnegative initial datum $u_0\in \LL^1_{\Phi_1}(\Omega)$.

Then $u\in C^{1,\alpha}_t((0,\infty)\times\Omega)$ for some $\alpha>0$ that depends only on $s$ and $m$.
Moreover, for any compact set $K\subset\subset\Omega$, and any $0<t_0<T$, we have
\begin{equation}\label{Thm6.reg}
\|u\|_{C^{1,\alpha}_t([t_0,T]\times K)}\leq C,
\end{equation}
where $C$ depends only on $N,s,m,\Omega, t_0,\|u_0\|_{\LL^1_{\Phi_1}(\Omega)}$, and $K$.
\end{thm}

\begin{rem}
A possible value for the exponent $\alpha$ in the previous theorem is $\alpha=\min\left\{\frac{1}{2m},1-s\right\}$, see Remark $\ref{rmk:alpha}$.
\end{rem}

Notice that results in Theorems \ref{Thm4}, \ref{Thm5}, and \ref{Thm6} imply that  \textit{solutions to \eqref{Dirichlet.Problem} are classical }  for any nonnegative initial datum $u_0\in \LL^1_{\Phi_1}(\Omega)$.

One may wonder whether the interior regularity in time can be improved to $C^\infty_t$.
This  is a very delicate problem that remains open, and is closely related to the higher order \emph{boundary} regularity in~$t$.
Indeed, recall that for nonlocal parabolic equations the interior regularity in the $t$-variable depends strongly on the regularity of the solution \emph{in all of} $\R^N$.
For example, for solutions to the fractional heat equation $u_t+(-\Delta)^su=0$ in $(0,1)\times B_1$, one has estimates of the form
\[\|u\|_{C^{k,\alpha}_x((\frac12,1)\times B_{1/2})}\leq C\|u\|_{L^\infty((0,1)\times\R^N)},\]
for all $k\geq0$ and $\alpha\in(0,1)$, which means that solutions are always $C^\infty$ in $x$.
However, analogous estimates in time do \emph{not} hold for $k\geq1$ and $\alpha\in(0,1)$.
Indeed, one can construct a solution to $u_t+(-\Delta)^su=0$ in $(0,1)\times B_1$, which is bounded in all of $\R^N$, but which is not $C^1$ in $t$ in $(\frac12,1)\times B_{1/2}$; see \cite[Section 2.4.1]{counterexample}.
Because of this delicate issue, since we only have a global H\"older regularity in the $t$-variable (by Theorem \ref{Thm4}) then we cannot prove more than $C^{1,\alpha}$ regularity in time in the interior of the domain.
We remark that, to our best knowledge, the higher boundary regularity in time is unknown even in the classical case $s=1$.

\section{Proof of the global lower bounds and Harnack inequalities}
\label{sect:lower}
In \cite{BV-PPR1,BV-PPR2-2} the authors prove that after some waiting time $t_*>0$ solutions are positive and satisfy sharp quantitative lower bounds. More precisely, letting
\begin{equation}\label{t*}
t_*=t_*(u_0)=\frac{c_*}{\|u_0\|_{\LL^1_{\Phi_1}(\Omega)}^{m-1}}\,,
\end{equation}
then the following global  lower bounds holds:
\begin{equation}\label{lower.bdds.BV-PPR1}
u(t,x)\ge c_0\frac{\Phi_1(x)^{\frac{1}{m}}}{t^{\frac{1}{m-1}}}\,,\qquad\mbox{for any $t\ge t_*$ and all $x\in \Omega$\,,}
\end{equation}
where the positive constants $c_*,c_0$ depend on only on $N,s,m,\Omega$. The above result holds true for any $0<s \le 1$. Recall that $\Phi_1\asymp d^s=\dist(\cdot,\partial\Omega)^s$, see \eqref{Phi.Estimates}.

It was posed in \cite{BV-PPR1} as an open problem to find precise lower bounds for small times, if possible. Here we solve that open problem.\\

\noindent\textbf{Infinite vs finite speed of propagation. }The purpose of this section is to extend the above lower bound to all $t\in (0,t_*)$.
As mentioned in the previous sections, this can be done only when $0<s<1$, and it appears in a relevant way the role of the nonlocal operator $\A$.

We show here, in a quantitative way, that when $s<1$ the fractional PME has infinite speed of propagation; this is obtained by the lower bounds of Theorem \ref{Thm1}, namely \begin{equation}\label{lower.thm.1.recall}
u(t,x)\ge c(u_0,t)\Phi_1(x)^{\frac{1}{m}}\,,\qquad\mbox{for any $0\le t\le t_*$ and all $x\in \Omega$\,.}
\end{equation}
Notice that, while for $t \geq t_*$ the constant in \eqref{lower.bdds.BV-PPR1} does not depend on the initial datum,
the constant appearing in the lower bound above does depend on $u_0$\,.

\noindent\textbf{Approximate solutions and their properties. }We begin by defining the following class of approximate solutions, that we will use throughout the proof.
Let us fix $\delta>0$ and consider a ``larger'' approximate problem:
\begin{equation}\label{probl.approx.soln}
\left\{
\begin{array}{lll}
\partial_t u_\delta=-\A u_\delta^m & \qquad\mbox{for any $(t,x)\in (0,\infty)\times \Omega$}\\
u_\delta(t,x)=\delta &\qquad\mbox{for any $(t,x)\in (0,\infty)\times (\RR^N\setminus\Omega)$}\\
u_\delta(0,x)=u_0(x)+\delta &\qquad\mbox{for any $x\in\Omega$}\,.
\end{array}
\right.
\end{equation}
We summarize here below the basic properties of $u_\delta$\,, which follows by adaptation of the theory in \cite{BV-PPR1, BV-PPR2-1, BV-PPR2-2, DPQRV2, VDPQR} that we recall in Appendix II for reader's convenience. We prefer to list only the properties of the approximate solutions that we need, in order not to break the flow of exposition.

Approximate solutions $u_\delta$ exists, are unique and bounded for all $(t,x)\in (0,\infty)\times\overline{\Omega}$\,, when $0\le u_0\in \LL^1_{\Phi_1}(\Omega)$\,.
Also, they are uniformly positive: for any $t \geq 0$,
\begin{equation}\label{approx.soln.positivity}
u_\delta(t,x)\ge \delta>0 \qquad\mbox{for a.e. $x \in \Omega$.}
\end{equation}
This implies that the equation for $u_\delta$ is never degenerate and solutions are smooth in the interior
(say, $C^\infty$ in space and $C^1$ in time). A possible way to prove this fact is by following the very same argument that we use in the next sections to prove the smoothness of $u$, cf. Section \ref{app.reg.approx} in Appendix II for some more details.

Also, by a comparison principle, for all $\delta>0$ the approximate solution $u_\delta$ are ordered and they lie above $u$: more precisely,
for all $\delta>\delta'>0$ and $t \geq 0$,
\begin{equation}\label{approx.soln.comparison1}
u_\delta(t,x)\ge u_{\delta'}(t,x) \qquad\mbox{for $x \in \Omega$}
\end{equation}
and
\begin{equation}\label{approx.soln.comparison}
u_\delta(t,x)\ge u(t,x) \qquad\mbox{for a.e. $x \in \Omega$\,.}
\end{equation}
Finally, as the following lemma shows, they converge in $\LL^1_{\Phi_1}(\Omega)$ to $u$, as $\delta \to 0$.

\begin{lem}\label{Lem.0}
Assume $0\le u_0\in\LL^1_{\Phi_1}(\Omega)$. Then, for any $t\geq 0$, we have
$$
\|u_\delta(t)-u(t)\|_{\LL^1_{\Phi_1}(\Omega)}\le \|u_\delta(0)-u_0\|_{\LL^1_{\Phi_1}(\Omega)}=\delta\,\|\Phi_1\|_{\LL^1(\Omega)}\,.
$$
\end{lem}
\noindent {\sl Proof.~}We provide a formal proof of this Lemma for reader's convenience. The proof can be made rigorous by taking smooth approximations of the test function $\Phi_1(x)\chi_{[0,t]}(\tau)$ and then passing to the limit, as it has been done, for instance, in the proof of \cite[Proposition 4.2]{BV-PPR1} or \cite[Proposition 5.1]{BV-PPR2-1}.

Since $\Phi_1$ is the first eigenfunction for $\A$ we know that $\A\Phi_1=\lambda_1\Phi_1$
for some $\lambda_1>0$ (see Appendix I).
Hence, using the  equation for both $u$ and $u_\delta$ and an integration by parts, we get
\[\begin{split}
\int_\Omega (u_\delta(t,x)-u(t,x))\Phi_1(x)\dx
&-\int_\Omega (u_\delta(0,x)-u_0(x))\Phi_1(x)\dx
=\int_0^t\int_\Omega\Phi_1\partial_t(u_\delta-u)\dx\dt\\
&=-\int_0^t\int_\Omega\Phi_1\, \A(u_\delta^m-u^m)\dx\dt
=-\int_0^t\int_\Omega\A\Phi_1\, (u_\delta^m-u^m)\dx\dt\\
&=-\lambda_1\int_0^t\int_\Omega \Phi_1\, (u_\delta^m-u^m)\dx\dt\le 0\\
\end{split}
\]
where at the last step we used \eqref{approx.soln.comparison}.
Recalling \eqref{approx.soln.comparison}, that $u_\delta(0,\cdot)-u_0(\cdot)=\delta$,
and that $\int_\Omega\Phi_1=1$, the result follows.\qed

Notice that, as a consequence of \eqref{approx.soln.comparison1} and Lemma \ref{Lem.0}, we deduce that $u_\delta$ converge pointwise to $u$
at almost every point: more precisely, for all $t\geq 0,$
\begin{equation}\label{limit.sol}
{u}(t,x)=\lim_{\delta\to 0^+}u_\delta(t,x)\ge 0\qquad\mbox{for a.e. $x \in \Omega$\,.}
\end{equation}
Next we prove a lower bound for $\LL^p_{\Phi_1}(\Omega)$ norms, a crucial ingredient of the proof of Theorem \ref{Thm1}, and that may have its own interest.
\begin{lem}\label{lem.abs.lower}
Let $u$ be a solution to problem \eqref{Dirichlet.Problem} corresponding to the initial datum $u_0\in\LL^1_{\Phi_1}(\Omega)$. Then the following lower bound holds true for any $t\in [0,t_*]$ and $p\ge 1$:
\begin{equation}\label{lem1.lower.bdd}
c_2 \left(\int_{\Omega}u_0(x)\Phi_1(x)\dx \right)^p \le \int_{\Omega}u^p(t,x)\Phi_1(x)\dx\le \int_{\Omega}u_\delta^p(t,x)\Phi_1(x)\dx.
\end{equation}
Here $t_*=c_*\|u_0\|_{\LL^1_{\Phi_1}(\Omega)}^{-(m-1)}$ is as in \eqref{t*},
and $c_2, c_*>0$ are positive constants that depend only on $N,s,m,p,\Omega$. Notice that $c_*$ has explicit form given in {\rm \cite{BV-PPR1,BV-PPR2-1, BV-PPR2-2}}, while the form of $c_2$ is given in the proof.
\end{lem}
\noindent {\sl Proof.~}Let us first recall the estimates of \cite[Corollary 6.2]{BV-PPR1} or \cite[Proposition 8.1]{BV-PPR2-1}. For any time $t_0$ satisfying
\begin{equation}\label{lem1.tsmall}
0\le t_0 \le c_1\|u_0\|_{\LL^1_{\Phi_1}(\Omega)}^{-(m-1)}=\frac{c_1}{c_*}t_*
\end{equation}
we have that
\begin{equation}\label{lem1.1}
\frac{1}{2}\int_{\Omega}u_0(x)\,\Phi_1(x)\dx\le \int_{\Omega}u(t_0,x)\,\Phi_1(x)\dx\,,
\end{equation}
where $c_1>0$ is a positive constant that depends only on $N,s,m,\Omega,$ and has explicit form given in \cite{BV-PPR1,BV-PPR2-1, BV-PPR2-2}.
Therefore, by H\"older inequality, the proof of \eqref{lem1.lower.bdd} is complete for all $t\in [0,t_0]$ provided $t_0$ satisfies \eqref{lem1.tsmall}.

We now prove the inequality for $t\in[t_0,t_*]$.
Recalling the monotonicity estimates of Benilan and Crandall \cite{BCr}, namely that $t^{\frac{1}{m-1}}u(t,x)\ge t_0^{\frac{1}{m-1}}u(t_0,x)$ for almost every $x\in\Omega$\,, we conclude that
\begin{equation}\label{lem1.2}\begin{split}
\frac{1}{2}\int_{\Omega}u_0(x)\,\Phi_1(x)\dx
    &\le \int_{\Omega}u(t_0,x)\,\Phi_1(x)\dx
     \le \left(\frac{t}{t_0}\right)^{\frac{1}{m-1}}\int_{\Omega}u(t,x)\,\Phi_1(x)\dx \\
    &\le \left(\frac{t}{t_0}\right)^{\frac{1}{m-1}}\|\Phi_1\|_{\LL^1(\Omega)}^{1-\frac{1}{p}}
        \left(\int_{\Omega}u^p(t,x)\,\Phi_1(x)\dx\right)^{\frac{1}{p}} \\
    &\le \left(\frac{t_*}{t_0}\right)^{\frac{1}{m-1}}\|\Phi_1\|_{\LL^1(\Omega)}^{1-\frac{1}{p}}
        \left(\int_{\Omega}u^p(t,x)\,\Phi_1(x)\dx\right)^{\frac{1}{p}} \\
    &= \left(\frac{2c_*}{c_1}\right)^{\frac{1}{m-1}}\|\Phi_1\|_{\LL^1(\Omega)}^{1-\frac{1}{p}}
        \left(\int_{\Omega}u^p(t,x)\,\Phi_1(x)\dx\right)^{\frac{1}{p}} \\
    &:=\frac{c_2^{\frac{1}{p}}}{2} \left(\int_{\Omega}u^p(t,x)\,\Phi_1(x)\dx\right)^{\frac{1}{p}}
\end{split}\end{equation}
where in the last step we have taken $t_0=c_1t_*/(2c_*)$\,, which is compatible with the restriction \eqref{lem1.tsmall}. This concludes the proof once we recall that $u_\delta\ge u$\,.\qed


\subsection{Proof of Theorem \ref{Thm1}}
We shall prove the lower bound \eqref{Thm1.lower.bdd} first for the approximate solutions $u_\delta$\,, i.e., solutions to problem \eqref{probl.approx.soln}  keeping $\delta>0$ fixed, and then we let $\delta\to 0^+$\,. We split the argument into several steps.

\noindent$\bullet~$\textsc{Step 1. }\textit{Reduction to an approximate problem. }Let us fix $\delta>0$ and consider the approximate solutions $u_\delta$. Recall that $u_\delta\ge \delta>0$ (see \eqref{approx.soln.positivity}), that $u_\delta \geq u$
(see \eqref{approx.soln.comparison}), and that $u$ coincides with their monotone limit in the strong $\LL^1_{\Phi_1}(\Omega)$ topology (cf. Lemma \ref{Lem.0}). We will keep $\delta>0$ fixed until the very last step.

\noindent$\bullet~$\textsc{Step 2. }\textit{Positivity for the approximate problem: Lower barriers. }We will show first that our positivity result holds for approximate solutions $u_\delta$, with a lower bound independent of $\delta$. Assume that $0\le t\le t_*$, and define the function (lower barrier)
\begin{equation}\label{lower.barrier}
\psi(t,x)=\kappa_0\,t\,\Phi_1(x)^{\frac{1}{m}}\,,
\end{equation}
where $\kappa_0>0$ is a positive parameter that will be fixed later (as a function of the initial data $u_0$). In order to determine $\kappa_0$, we use the lower bounds of \cite{BV-PPR1, BV-PPR2-1, BV-PPR2-2}, see also \eqref{lower.bdds.BV-PPR1}: for all $t\ge t_*$ we have $u(t,x)\ge c_0 \Phi_1(x)^{\frac{1}{m}} t^{-\frac{1}{m-1}}\,,$ where $t_*=t_*(u_0)>0$ is given in \eqref{t*}. Therefore\,, since we already know that $u_\delta\ge u$\,, we have
\begin{equation}\label{ineq.t*}
\psi(t_*,x)=\kappa_0\,t_*\,\Phi_1(x)^{\frac{1}{m}}< c_0\frac{\Phi_1(x)^{\frac{1}{m}}}{t^{\frac{1}{m-1}}}\le u(t_*,x)\le u_{\delta}(t_*,x)
\end{equation}
whenever
\begin{equation}\label{cond.kappa0}
0<\kappa_0< \frac{c_0}{t_*^{\frac{m}{m-1}}}=c_0 c_*^{-\frac{m}{m-1}}\,\|u_0\|_{\LL^1_{\Phi_1}(\Omega)}^{m}=\overline{\kappa}_0\,.
\end{equation}
\textsl{Claim. }We want to prove that $\psi(t,x)< u_{\delta}(t,x)$ for all $0\le t\le t_*$\,, all $x\in \Omega$, for a sufficiently small  $\kappa_0$ independent of $\delta$ that satisfies also the condition \eqref{cond.kappa0}\,.

We will prove the claim by contradiction. Since $\psi(0)=0<\delta\leq u_\delta(0)$,
Assume that the inequality $\psi<u_\delta$ is false in $[0,t_*]\times \overline\Omega$
 and let $(t_\delta,x_\delta)$ be the first contact point.
 Notice that  $0<t_\delta<t_*$ (since $\psi(0,x)=0<\delta\leq u_\delta(0,x)$ and by  \eqref{ineq.t*}).
 Also, $(t_\delta,x_\delta)$ cannot belong to $(0,t_*)\times\partial\Omega$ since $u_\delta(t,x)=\delta>0=\psi(t,x)$ there,
thus $(t_\delta,x_\delta) \in (0,t_*)\times \Omega$.
Now, since $(t_\delta,x_\delta)$ is the first contact point, we necessarily have that
\begin{equation}\label{contact.1}
u_\delta(t_\delta,x_\delta)= \psi(t_\delta,x_\delta)\qquad\mbox{and}\qquad
u_\delta(t,x)\ge \psi(t,x)\,\quad\forall t\in [0,t_\delta]\,,\;\; \forall x\in\Omega\,,
\end{equation}
and, as a consequence,
\begin{equation}\label{contact.2}
\partial_t u_\delta(t_\delta,x_\delta)\le \partial_t \psi(t_\delta,x_\delta)=\kappa_0\,\Phi_1(x)^{\frac{1}{m}}\,.
\end{equation}
Thus, using  \eqref{contact.2} and
\[
\A(\psi^m)(t,x)= \kappa_0^m\,t^m\, \A(\Phi_1)(x)= \kappa_0^m\,t^m\, \lambda_1\,\Phi_1(x)\,,
\]
 we first establish an upper bound for $-\A(u_\delta^m-\psi^m)(t_\delta,x_\delta)$ as follows:
\begin{equation}\label{contact.3}
-\A(u_\delta^m-\psi^m)(t_\delta,x_\delta)=\partial_t u(t_\delta,x_\delta)+\A(\psi^m)(t_\delta,x_\delta)
\le \kappa_0\,\Phi_1(x_\delta)^{\frac{1}{m}}(1+\lambda_1\, \kappa_0^{m-1} t_\delta^m ).
\end{equation}
Next, we want to prove lower bounds for $-\A(u_\delta^m-\psi^m)(t_\delta,x_\delta)$
(this is the crucial point where the nonlocality of the operator enters).
To this end we recall that\\
\begin{equation}\label{Operator.Hyp.lower}
\A(f)(x)=\int_{\RR^N}\big(f(x)-f(y)\big)\K(x,y)\dy\qquad\mbox{and}\qquad \inf_{x,y\in \Omega}\K(x,y)\ge \underline{\kappa}_\Omega>0\,.
\end{equation}
Hence, using  \eqref{contact.1} and \eqref{Operator.Hyp.lower}, we get
\begin{equation}\label{contact.4}\begin{split}
-\A(u_\delta^m-\psi^m)(t_\delta,x_\delta)
&=-\int_{\RR^N}\left[(u_\delta^m(t_\delta,x_\delta)-u_\delta^m(t_\delta,y))-(\psi^m(t_\delta,x_\delta)-\psi^m(t_\delta,y))\right]\K(x,y)\dy\\
&=\int_{\Omega}\left[u_\delta^m(t_\delta,y)-\psi^m(t_\delta,y)\right]\K(x,y)\dy\\
&\ge \underline{\kappa}_\Omega\,\int_{\Omega}\left[u_\delta^m(t_\delta,y)-\psi^m(t_\delta,y)\right]\dy\\
&= \underline{\kappa}_\Omega\,\int_{\Omega}u_\delta^m(t_\delta,y)\dy-\underline{\kappa}_\Omega\,\int_{\Omega}\psi^m(t_\delta,y)\dy.
\end{split}
\end{equation}
Combining the upper and lower bounds \eqref{contact.3} and \eqref{contact.4} we obtain
\begin{equation}\label{contact.5}\begin{split}
\underline{\kappa}_\Omega\,\int_{\Omega}u_\delta^m(t_\delta,y)\dy
& \le \kappa_0\,\Phi_1(x_\delta)^{\frac{1}{m}}(1+\lambda_1 \kappa_0^{m-1}\,t_\delta^m )
    +\underline{\kappa}_\Omega\,\int_{\Omega}\psi^m(t_\delta,y)\dy\\
& = \kappa_0\,\Phi_1(x_\delta)^{\frac{1}{m}}(1+\lambda_1 \kappa_0^{m-1}\,t_\delta^m )+\underline{\kappa}_\Omega\,\kappa_0^m t_\delta^m \|\Phi_1\|_{\LL^1(\Omega)}\\
& = \kappa_0\,\left[\Phi_1(x_\delta)^{\frac{1}{m}}(1+\lambda_1 \kappa_0^{m-1}\,t_\delta^m )
        +\underline{\kappa}_\Omega\,\kappa_0^{m-1} t_\delta^m\|\Phi_1\|_{\LL^1(\Omega)}\right]\\
& = \kappa_0\,\left[\|\Phi_1\|_{\LL^\infty(\Omega)}^{\frac{1}{m}}(1+\lambda_1\, \overline{\kappa}_0^{m-1} t_*^m )
        +\underline{\kappa}_\Omega\,\overline{\kappa}_0^{m-1} t_*^m\|\Phi_1\|_{\LL^1(\Omega)}\right]
        \le \underline{\kappa}_\Omega\, \overline{K} \,\kappa_0\,.\\
\end{split}\end{equation}
 In the very last step we have used that $\overline{\kappa}_0^{m-1} t_*^m=c_0^{m-1}$ and that $\|\Phi_1\|_{\LL^1(\Omega)}\le \|\Phi_1\|_{\LL^\infty(\Omega)}|\Omega|$, to get
\begin{equation}\label{kappa01}\begin{split}
\|\Phi_1\|_{\LL^\infty(\Omega)}^{\frac{1}{m}} &(1+\lambda_1\, \overline{\kappa}_0^{m-1} t_*^m )
        +\underline{\kappa}_\Omega\,\overline{\kappa}_0^{m-1} t_*^m\|\Phi_1\|_{\LL^1(\Omega)}\\
&\le \|\Phi_1\|_{\LL^\infty(\Omega)}^{\frac{1}{m}}(1+\lambda_1\,c_0^{m-1})
        +\underline{\kappa}_\Omega\,c_0^{m-1}\|\Phi_1\|_{\LL^1(\Omega)}:=\overline{K}\,.
\end{split}\end{equation}
As a consequence, for $t_\delta\in (0,t_*)$ we have
\begin{equation}\label{contact.6}
\int_{\Omega}u_\delta^m(t_\delta,y)\dy\le  \overline{K}  \,\kappa_0\,,\qquad\mbox{for all $0< \kappa_0<\overline{\kappa}_0$\,.}
\end{equation}
To get a contradiction we need to use the lower bounds of Lemma \ref{lem.abs.lower}: choosing $p=m>1$ we get that, for all $0\le t_\delta\le t_*$,
\begin{equation}\label{lem1.lowe.bdds.1}
c_2 \left(\int_{\Omega}u_0(x)\Phi_1(x)\dx \right)^m \le   \int_{\Omega}u_\delta^m(t_\delta,x)\Phi_1(x)\dx
\le \|\Phi_1\|_{\LL^\infty(\Omega)} \int_{\Omega}u_\delta^m(t_\delta,x)\dx
\end{equation}
where $c_2$ is an explicit positive constant that depends only on $N,s,m,\Omega$.

Combining now \eqref{contact.6} and \eqref{lem1.lowe.bdds.1} we obtain
$$
\frac{c_2}{\|\Phi_1\|_{\LL^\infty(\Omega)}} \left(\int_{\Omega}u_0(x)\Phi_1(x)\dx \right)^m\le  \overline{K} \,\kappa_0\,,
\qquad\mbox{for all $0< \kappa_0<\overline{\kappa}_0$\,.}
$$
hence a contradiction (choosing $\kappa_0$ small enough) whenever $u_0\not\equiv 0$ a.e. in $\Omega$\,, independently of $\delta>0$. More precisely we can choose $\kappa_0$ to be
\begin{equation}\label{kappa0}
\kappa_0< \min\left\{\frac{c_2\,  \overline{K}^{-1}  }{\|\Phi_1\|_{\LL^\infty}}\,, \frac{c_0}{c_*^{\frac{m}{m-1}}}\right\}\|u_0\|_{\LL^1_{\Phi_1}}^m,
\end{equation}
which proves the claim.   The expression of $\overline{K}$ is given in \eqref{kappa01}.

Summing up, we have proved the claim, in the following precise form: for any $\delta>0$
\begin{equation}\label{step.2.final}
u_\delta(t,x)\ge \psi(t,x)=\kappa_0\,t\,\Phi_1(x)^{\frac{1}{m}}\qquad\mbox{for all $0<t\le t_*$ and $x\in \Omega$\,,}
\end{equation}
whenever
\begin{equation}\label{step.2.final.cond}
0<\kappa_0< \min\left\{\frac{c_2\,   \overline{K}^{-1}  }{\|\Phi_1\|_{\LL^\infty}}\,, \frac{c_0}{c_*^{\frac{m}{m-1}}}\right\}\|u_0\|_{\LL^1_{\Phi_1}}^m, \qquad\mbox{with $\overline{K}$ given in \eqref{kappa01}. }
\end{equation}

\noindent$\bullet~$\textsc{Step 3. }\textit{Taking the limit $\delta\to 0^+$: Quantitative positivity for solutions.}
Using \eqref{limit.sol} and letting $\delta\to 0^+$ in \eqref{step.2.final}, we obtain that \eqref{Thm1.lower.bdd} follows with
\begin{equation}\label{kappa}
\kappa:=\frac{\kappa_0}{\|u_0\|_{\LL^1_{\Phi_1}}^m}=  \min\left\{\frac{c_2\,\overline{K}_\Omega}{\|\Phi_1\|_{\LL^\infty(\Omega)}}\,,
\frac{1}{\overline{K}_\Omega}\frac{c_2}{\|\Phi_1\|_{\LL^\infty(\Omega)}}\,,
\frac{c_0}{c_*^{\frac{m}{m-1}}}\right\}\,.\mbox{\qed}
\end{equation}

\subsection{Proofs of Theorems \ref{Thm2} and \ref{Thm3}}

\noindent\textbf{Proof of Theorem \ref{Thm2}. }When $0\le t\le t_*$ we have that inequality \eqref{Thm1.lower.bdd} can be rewritten as follows
\begin{equation}\label{lower.bdd.thm2}
u(t,x)\ge \kappa\|u_0\|^m_{\LL^1_{\Phi_1}(\Omega)}\,t\,\Phi_1(x)^{\frac{1}{m}}
=\kappa\|u_0\|^m_{\LL^1_{\Phi_1}(\Omega)}\,t^{\frac{m}{m-1}}\,\frac{\Phi_1(x)^{\frac{1}{m}}}{t^{\frac{1}{m-1}}}
=\frac{\kappa}{c_*^{\frac{m}{m-1}}} \frac{t^{\frac{m}{m-1}}}{t_*^{\frac{m}{m-1}}}\frac{\Phi_1(x)^{\frac{1}{m}}}{t^{\frac{1}{m-1}}},
\end{equation}
where $t_*=c_*\|u_0\|_{\LL^1_{\Phi_1}(\Omega)}^{-(m-1)}$.
We combine this lower bound with the one valid for large times, namely inequality \eqref{thm.GHP.PME.1} that holds for any $t\ge t_*$ and all $x\in \overline{\Omega}$\,, and reads
\begin{equation}\label{thm.GHP.PME.1.thm2}
u(t,x)\ge H_0\,\frac{\Phi_1(x)^{\frac{1}{m}}}{t^{\frac{1}{m-1}}}\,.
\end{equation}
Joining the two lower bounds \eqref{lower.bdd.thm2} and \eqref{thm.GHP.PME.1.thm2} finally gives
\[
u(t,x)\ge \min\left\{H_0, \frac{\kappa}{c_*^{\frac{m}{m-1}}} \frac{t^{\frac{m}{m-1}}}{t_*^{\frac{m}{m-1}}}\right\}\frac{\Phi_1(x)^{\frac{1}{m}}}{t^{\frac{1}{m-1}}}
:=\underline{\kappa}\left(1\wedge\frac{t}{t_*}\right)^{\frac{m}{m-1}} \frac{\Phi_1(x)^{\frac{1}{m}}}{t^{\frac{1}{m-1}}}\,,
\]
which is exactly \eqref{Thm2.GHP}.\qed

\noindent\textbf{Proof of Theorem \ref{Thm3}. }
We recall the Harnack inequality for the first eigenfunction
\begin{equation}\label{Harnack.Phi1}
\sup_{x\in B_R(x_0)} \Phi_1(x)\le H_{R,x_0} \inf_{x\in B_R(x_0)} \Phi_1(x)\,,\qquad H_{R,x_0}=H_{N,s,\Omega}\,\dist(B_R(x_0),\partial\Omega)^{-s},
\end{equation}
which follows from \eqref{Phi.Estimates}. We then use the Global Harnack Principle, namely inequality \eqref{Thm2.GHP} for all $t>0$, as follows:
\[
\sup_{x\in B_R(x_0)} u(t,x)\le  \frac{\overline{\kappa}}{t^{\frac{1}{m-1}}}  \sup_{x\in B_R(x_0)}\Phi_1(x)^{\frac{1}{m}}
\le \frac{\overline{\kappa}\,H_{R,x_0}}{t^{\frac{1}{m-1}}}  \inf_{x\in B_R(x_0)}\Phi_1(x)^{\frac{1}{m}}
\le \frac{\overline{\kappa}\,H_{R,x_0}\underline{\kappa}}{\left(1\wedge\frac{t}{t_*}\right)^{\frac{m}{m-1}}}\inf_{x\in B_R(x_0)} u(t,x)\,,
\]
where in the second inequality we have used \eqref{Harnack.Phi1}. This is exactly \eqref{Thm3.local.Harnack} with $\mathcal{H}=\overline{\kappa}\,H_{R,x_0}\underline{\kappa}$\,.\qed

\section{Proofs of the regularity results}

\noindent\textit{About the class of solutions considered. }In this section we will work with weak solutions, a class of solution contained in the class of weak dual solutions to which the results of Theorems \ref{Thm1}\,, \ref{Thm2} and \ref{Thm3} apply.
On the other hand, arguing by approximation, all results that we shall prove for weak solutions will also be valid for weak dual solution.
To clarify this point, we recall that in \cite{BV-PPR1} the authors construct weak dual solutions starting from weak solutions when the operator $\A$ is the so-called Spectral Fractional Laplacian (which is different from the operators that we use here), but the methods can be adapted to this setting as well\,, cf. also \cite{BSV2013}.
More precisely, as explained in \cite{BV-PPR1}, weak dual solutions have to be interpreted as a class of limit solutions which satisfy a quasi $\LL^1_{\Phi_1}$ contraction, cf. \cite{BV-PPR2-1}\,; as a consequence, the minimal weak dual solution obtained by approximation from below turns out to be unique as an $\LL^1_{\Phi_1}$ function. In this section we will prove our regularity results for weak solutions. To conclude that the same regularity results for weak dual solutions, it  suffices to approximate a weak dual solution from below with a sequence of weak solutions to which the regularity results apply, and notice that all the constants in the regularity estimates are stable under this limit process.
Indeed, these estimates only depend on the supremum of $u$, and the latter is controlled by the $\LL^1_{\Phi_1}$-norm of the initial datum in view of the smoothing effects, cf. \cite{BV-PPR1,BV-PPR2-1}. Hence the constants are stable under this approximation process and the regularity results valid for weak solutions can be extended to weak dual solutions.


In this section we prove first Theorem \ref{Thm4} and then Theorem \ref{Thm5}.
The main step towards Theorem \ref{Thm4} is the following regularity estimate.

\begin{prop}\label{prop.regularity}
Let $w$ be a weak solution of
\begin{equation}\label{eq-regularity}
w_t+(-\Delta)^sw^m=0\quad for\ (t,x)\in(0,1)\times B_1.
\end{equation}
Assume that for all $t\in(0,1)$ we have
\begin{equation}\label{growth}
|w(t,x)|\leq C_0\bigl(1+|x|^{s/m}\bigr)\quad for\ x\in\RR^N,\qquad\textrm{and}\qquad w(t,x)\geq\delta>0\quad  \textrm{for}   \ x\in B_3\,,
\end{equation}
for some constants $C_0$ and $\delta>0$.
Then, for each $0<\beta\leq s$ we have the estimate
\[\|w\|_{C^{\beta,\beta/2s}_{x,t}((\frac12,1)\times B_{1/2})}\leq C,\]
for some constant $C$ that depends only on $N$, $s$, $m$, $\beta$, $\delta$, and $C_0$.
\end{prop}

We will split the proof of Proposition \ref{prop.regularity} into some steps.
First, we show the following.

\begin{lem}\label{lem.regularity1}
Let $w$ be a weak solution of \eqref{eq-regularity} satisfying \eqref{growth}.
Then, for some small $\bar\alpha>0$ we have the estimate
\[\|w\|_{C^{\bar\alpha,\bar\alpha/2s}_{x,t}((\frac12,1)\times B_{1/2})}\leq C,\]
for some $C$ that depends only on $N$, $s$, $m$, $\delta$, and $C_0$.
\end{lem}

As we will see, the previous Lemma is a consequence of the following H\"older estimate by Felsinger-Kassmann \cite{FK}.

\begin{thm}[\cite{FK}]\label{FK-estimate}
Let $L$ be an integro-differential operator defined by
\begin{equation}\label{op-kassmann}
Lw(t,x)=PV\int_{\RR^N}\bigl(w(t,x)-w(t,y)\bigr)\frac{a(t,x,y)}{|x-y|^{n+2s}}\,dy,
\end{equation}
where $a$ satisfies $a(t,x,y)=a(t,y,x)$ and $\Lambda^{-1}\leq a(t,x,y)\leq \Lambda$ for some $\Lambda\geq1$.
Let $w\in L^\infty((0,1)\times\RR^N)$ be any weak solution of
\[w_t+Lw=f\qquad \text{for}\ (t,x)\in(0,1)\times B_1,\]
with $f\in L^\infty((0,1)\times B_1)$.
Then, for some $\bar\alpha>0$ small we have the estimate
\[\|w\|_{C^{\bar\alpha,\bar\alpha/2s}_{x,t}((\frac12,1)\times B_{1/2})}\leq C\bigl(\|f\|_{L^\infty((0,1)\times B_1)}+\|w\|_{L^\infty((0,1)\times\RR^N)}\bigr).\]
The constants $C$ and $\bar\alpha$ depend only on $N$, $s$, and $\Lambda$.
\end{thm}

\begin{rem}
In {\rm \cite{FK}} the H\"older regularity is stated for simplicity for the case $f=0$.
Still, the weak Harnack inequality is established for the general case $f\in L^\infty$, and thus the proof of the H\"older regularity with a nonzero right hand side $f$ (as stated above) is just a minor modification of the one in~{\rm \cite{FK}}.
\end{rem}

Using this estimate, we can now prove Lemma \ref{lem.regularity1}.

\noindent {\bf Proof of Lemma \ref{lem.regularity1}.~}
First, let $\rho\in C^\infty_c(B_4)$ be a smooth cutoff function with $\rho\equiv1$ in $B_3$, and set $v:=w\rho$.
Then, we have
\[(-\Delta)^sv^m=(-\Delta)^sw^m-(-\Delta)^s\bigl((1-\rho^m)w^m\bigr).\]
Moreover, since $(1-\rho^m)w^m\equiv0$ in $B_3$,
\[-(-\Delta)^s\bigl((1-\rho^m)w^m\bigr)(x)=c_{n,s}\int_{\RR^N}(1-\rho^m(y))\,w^m(t,y)\,\frac{dy}{|x-y|^{n+2s}}\qquad\textrm{for}\quad x\in B_3\]
and this term belongs to $C^\infty(B_2)$, thus we have
\[v_t+(-\Delta)^sv^m=g(t,x)\in C^\infty(B_2).\]

Consider now $\eta\in C^\infty_c(B_2)$ a radial smooth cutoff function such that $\eta\equiv1$ in $B_1$.
Then, for $x\in B_1$ and $t\in(0,1)$ we may write
\begin{equation}\label{from-w^m-to-w}
v^m(t,x)-v^m(t,y)=\bigl(v(t,x)-v(t,y)\bigr)a(t,x,y)+h(t,x,y),
\end{equation}
with
\[a(t,x,y)=\frac{v^m(t,x)-v^m(t,y)}{v(t,x)-v(t,y)}\,\eta(x-y)+\bigl(1-\eta(x-y)\bigr)\]
and
\[h(t,x,y)=\bigl(v^m(t,x)-v^m(t,y)-v(t,x)+v(t,y)\bigr)\bigl(1-\eta(x-y)\bigr).\]
Notice that
\begin{equation}\label{a-alternative}
a(t,x,y)=m\,\eta(x-y)\int_0^1 \left[v(t,x)+\lambda \big((v(t,y)-v(t,x)\big)\right]^{m-1}d\lambda
+\bigl(1-\eta(x-y)\bigr),
\end{equation}
hence, using \eqref{growth}, we find that
\[\Lambda^{-1}\leq a(t,x,y)\leq \Lambda\quad \textrm{for}\quad x\in B_1\]
for some constant $\Lambda\geq1$ depending on $\delta$ and $C_0$.
Thus, using \eqref{from-w^m-to-w}, we see that $v$ is a weak solution of
\[v_t+Lv=g+f\quad\textrm{in}\ (0,1)\times B_1,\]
with $L$ as in \eqref{op-kassmann}, $g$ as above, and
\[f(t,x)=c_{n,s}\int_{\RR^N} h(t,x,y) \frac{dy}{|x-y|^{n+2s}}.\]
Moreover, using that $v$ is bounded and $v\geq\delta>0$ in $B_3$, we also find that
\begin{equation}\label{growth-h}
|h(t,x,y)|\leq C\quad\textrm{for}\quad x\in B_1,\ y\in \RR^N,\qquad\textrm{and}\qquad  h\equiv0\quad\textrm{for}\quad |x-y|\leq 1.
\end{equation}
This yields $\|f\|_{L^\infty((0,1)\times B_1)}\leq C$, and therefore, by Theorem \ref{FK-estimate},
\[\|v\|_{C^{\bar\alpha,\bar\alpha/2s}_{x,t}((1/2,1)\times B_{1/2})}\leq C\]
for some constants $\bar\alpha>0$ and $C$ that depend only on $n$, $s$, $m$, $\delta$, and $C_0$.
Since $v=w$ in $(0,1)\times B_3$, the same estimate holds for $w$, and thus the Lemma is proved.
\qed

Once we know that the solution $u$ is strictly positive and H\"older continuous in a ball, one may apply the arguments of \cite{VDPQR} to obtain higher regularity of~$u$.
Indeed, the results of \cite{VDPQR} are for global solutions in $\R^N$, but all the proofs therein have a local nature -- this is in fact explicitly stated in \cite[Section 7]{VDPQR}.
Thanks to this, we may use the same arguments to obtain higher regularity of solutions $u\in C^\alpha$ satisfying $0<\delta\leq u\leq C_0$ in a ball: more precisely, we have the following.

\begin{lem}\label{lem.regularity2}
Let $w$ be a weak solution of \eqref{eq-regularity} satisfying \eqref{growth}.
Assume in addition that $w\in C^{\gamma,\gamma/2s}_{x,t}((0,2)\times B_4)$ for some $\gamma\in [\gamma_0, s]$.
Then, for some small $\epsilon>0$ we have the estimate
\[\|w\|_{C^{\gamma+\epsilon,(\gamma+\epsilon)/2s}_{x,t}((1,2)\times B_2)}\leq C.\]
Here $C$ depends only on $N$, $s$, $m$, $\delta$, $C_0$, $\|w\|_{C^{\gamma,\gamma/2s}_{x,t}((0,2)\times B_4)}$,
and $\epsilon$ only on $\gamma_0$, $s$, $m$.
\end{lem}

\noindent {\bf Proof of Lemma \ref{lem.regularity2}.~}
As in \cite{VDPQR} we note that, given $(x_0,t_0)\in (1,2)\times B_2$,
after the time rescaling $t\mapsto \frac{t}{m\,w^{m-1}(x_0,t_0)}$
we have that $w$ solves the equation
\[w_t+(-\Delta)^s w= (-\Delta)^sf,\]
where
\[f(x,t):=w(x,t)-\frac{w^m(x,t)}{m\,w^{m-1}(x_0,t_0)}.\]
Then, denoting $Y=(x,t)$, $Y_0=(x_0,t_0)$, and $|Y|_s=(|x|^2+|t|^{1/s})^{1/2}$, since $w\in C^{\gamma,\gamma/2s}_{x,t}((0,2)\times B_4)$, it follows by \cite[Equation (4.1)]{VDPQR} that
\[|f(Y_1)-f(Y_2)|\leq C|Y_1-Y_2|_s^\gamma\,\max\bigl\{|Y_1-Y_0|^\epsilon,\,|Y_2-Y_0|^\epsilon\bigr\}\]
for all $Y_1$ and $Y_2$ in $(0,1)\times B_4$.

Now, we define $F=f\eta$, where $\eta\in C^\infty_c(B_4)$ is a cutoff function such that $\eta\equiv1$ in $B_3$, and
consider the unique solution $W$ to
\[W_t+(-\Delta)^sW=(-\Delta)^sF\qquad \textrm{in}\quad (0,\infty)\times\R^N,\]
with $W=0$ for $t=0$.
Then, by the results in \cite[Section 4]{VDPQR}, we have that
\[|W(Y_0+Y)+W(Y_0-Y)-2W(Y_0)|\leq C|Y|_s^{\gamma+\epsilon}\]
for all $Y$ such that $Y+Y_0,Y-Y_0\in (1,2)\times B_2$.

On the other hand, we notice that $v=w-W$ satisfies
\[v_t+(-\Delta)^sv=(-\Delta)^s(f-F)\qquad \textrm{in}\quad (0,2)\times B_3,\]
and $f-F\equiv0$ in $(0,2)\times B_3$.
Thus, the right hand side $(-\Delta)^s(f-F)$ is of class $C^\infty$ inside $(0,2)\times B_3$, and
by parabolic regularity for the fractional heat equation (see for example \cite{JX} and \cite{CKriv})
we deduce that $v\in C^{\infty,1}_{x,t}((1,2)\times B_2)$ (however, as mentioned in the introduction, we cannot deduce from this that $v$ is of class $C^\infty$ in $t$!).
Recalling that $w=v+W$, this yields
\[|w(Y_0+Y)+w(Y_0-Y)-2w(Y_0)|\leq C|Y|_s^{\gamma+\epsilon}\]
for some $\epsilon>0$.

Finally, since this can be done for all $(x_0,t_0)\in (1,2)\times B_2$, we find
\[\|w\|_{C^{\gamma+\epsilon,(\gamma+\epsilon)/2s}_{x,t}((1,2)\times B_1)}\leq C,\]
and thus the lemma is proved.
\qed

We can now prove Proposition \ref{prop.regularity}.

\noindent {\bf Proof of Proposition \ref{prop.regularity}.~}
First, by Lemma \ref{lem.regularity1} we have that
\[\|w\|_{C^{\bar\alpha,\bar\alpha/2s}_{x,t}((1/2,1)\times B_{1/2})}\leq C\]
for some $\bar\alpha>0$.
Then, applying Lemma \ref{lem.regularity2} we find
\[\|w\|_{C^{\bar\alpha+\epsilon,(\bar\alpha+\epsilon)/2s}_{x,t}((3/4,1)\times B_{1/4})}\leq C\]
for some fixed $\epsilon>0$.
Iterating Lemma \ref{lem.regularity2} finitely many time, we find
\[\|w\|_{C^{\beta,\beta/2s}_{x,t}((1-2^q,1)\times B_{1/2^q})}\leq C,\]
for some integer $q\geq2$.
By a standard covering argument this yields
\[\|w\|_{C^{\beta,\beta/2s}_{x,t}((1/2,1)\times B_{1/2})}\leq C,\]
thus the Proposition is proved.
\qed

We now show Theorem \ref{Thm4}.

\noindent {\bf Proof of Theorem \ref{Thm4}.~}
Let us fix a ball $B_r(x_0)$ such that $2r=\textrm{dist}(x,\partial\Omega)$, and define
\[u_r(t,x)=r^{-s/m}\,u\left(t_0+r^{s(1+\frac1m)}t,\,x_0+rx\right).\]
Then, it follows by Theorem \ref{Thm2} that
\[0<\delta\leq u_r(t,x)\leq C_0,\]
for all $t\in[0,1]$ and $x\in B_1$, with constants $\delta>0$ and $C_0$ independent of $r$ and $x_0$.
Furthermore, again by Theorem \ref{Thm2} we have
\[
u_r(t,x)\leq Cr^{-s/m} d(x_0+rx)^{s/m}
\]
for all $x\in \RR^N$ and all $t\geq t_0$, so that we get
\[u_r(t,x)\leq Cr^{-s/m}\textrm{dist}(x_0+rx,\partial\Omega)^{s/m}\leq C\bigl(1+|x|^{s/m}\bigr)\]
for all $x\in \R^N$, $t\geq t_0$.
Noticing that
\[\partial_t u_r+(-\Delta)^su_r^m=0\qquad \textrm{in}\quad (0,1)\times B_1,\]
we see that $u_r$ satisfies the hypotheses of Proposition \ref{prop.regularity}, and thus for any $\beta\leq s$ there exists a constant $C$ such that
\[\|u_r\|_{C^{\beta,\beta/2s}_{x,t}((1/2,1)\times B_{1/2})}\leq C.\]
Rescaling back to $u$, we find
\[[u]_{C^{\beta,\beta/2s}_{x,t}\left(\left(t_0+\frac12r^{s(1+\frac1m)},t_0+r^{s(1+\frac1m)}\right)\times B_{r/2}\right)}\leq Cr^{s/m-\beta}.\]
In particular, setting $\beta=s/m$, we find that, for all such ball $B_r(x_0)\subset B_{2r}(x_0)\subset \Omega$,
\[[u]_{C^{s/m,1/2m}_{x,t}\left(\left(t_0+\frac12r^{s(1+\frac1m)},t_0+r^{s(1+\frac1m)}\right)\times B_{r/2}\right)}\leq C\]
for some constant $C$ independent of $r$ and $x_0$.
This means that (see for example \cite{RS-Dir})
\[\|u\|_{C^{s/m,1/2m}_{x,t}([t_0,T]\times\overline\Omega)}\leq C,\]
and thus the Theorem is proved.
\qed

In order to show Theorem \ref{Thm5} we will need the following.

\begin{prop}\label{prop.regularityB}
Let $w$ be a weak solution of
\eqref{eq-regularity} satisfying \eqref{growth}.
Then, for each $\beta>0$ we have the estimate
\[\|w\|_{C^{\beta}_x((1/2,1)\times B_{1/2})}\leq C,\]
for some constant $C$ that depends only on $N$, $s$, $m$, $\beta$, $\delta$, and $C_0$.
\end{prop}

To prove this result, we first show the following lemma.

\begin{lem}\label{lem.regularity2B}
Let $w$ be any weak solution of \eqref{eq-regularity} satisfying \eqref{growth}.
Assume in addition that $w\in C^{\gamma}_{x}((0,2)\times B_4)$ for some $\gamma\in [\gamma_0,1]$ with $\gamma_0>0$.
Then, for some small $\epsilon>0$ we have the estimate
\[\|w\|_{C^{\gamma+\epsilon}_{x}((1,2)\times B_2)}\leq C.\]
Here $C$ depends only on $\gamma_0$, $N$, $s$, $m$, $\delta$, $C_0$, $\|w\|_{C^{\gamma}_{x}((0,2)\times B_4)}$,
and $\epsilon$ depends only on $\gamma_0$, $s$, $m$.
\end{lem}

\noindent {\bf Proof of Lemma \ref{lem.regularity2B}.~}
As in the proof of Lemma \ref{lem.regularity2}, given $(x_0,t_0)\in (1,2)\times B_2$,
after the time rescaling $t\mapsto \frac{t}{m\,w^{m-1}(x_0,t_0)}$
we have that $w$ solves the equation
\[w_t+(-\Delta)^s w= (-\Delta)^sf,\]
where
\[f(t,x):=w(t,x)-\frac{w^m(t,x)}{m\,w^{m-1}(t_0,x_0)},\]
and that
\[|f(t_0,x_1)-f(t_0,x_2)|\leq C|x_1-x_2|^\gamma\,\max\bigl\{|x_1-x_0|^\epsilon,\,|x_2-x_0|^\epsilon\bigr\}\]
for all $(t_0,x_1)$ and $(t_0,x_2)$ in $(0,1)\times B_4$.

Now, again following the argument in the proof of Lemma \ref{lem.regularity2}, we define $F=f\eta$, where $\eta\in C^\infty_c(B_4)$ is such that $\eta\equiv1$ in $B_3$, and we consider the unique solution $W$ to
\[W_t+(-\Delta)^sW=(-\Delta)^sF\qquad \textrm{in}\quad (0,\infty)\times\R^N,\]
with $W=0$ for $t=0$.
Then, by the results of \cite{VDPQR} (see the proofs \cite[Sections 4 and 5, and Proposition 6.1]{VDPQR}) we have that
\[|W(t_0,x_0+x)+W(t_0,x_0-x)-2W(t_0,x_0)|\leq C|x|^{\gamma+\epsilon}\]
for all $x$ such that $x_0\pm x\in B_2$.

On the other hand, we have that $v=w-W$ satisfies
\[v_t+(-\Delta)^sv=(-\Delta)^s(f-F)\qquad \textrm{in}\quad (0,2)\times B_3,\]
and $f-F\equiv0$ in $(0,2)\times B_3$.
Thus, the right hand side $(-\Delta)^s(f-F)$ belongs to $C^\infty_x$, and thus $v\in C^{\infty}_{x}((1,2)\times B_2)$.
Recalling that $w=v+W$, this yields
\[|w(t_0,x_0+x)+w(t_0,x_0-x)-2w(t_0,x_0)|\leq C|x|^{\gamma+\epsilon}\]
for some $\epsilon>0$.

Finally, since this argument can be  done at all points $(t_0,x_0)\in (1,2)\times B_2$, we find
\[\|w\|_{C^{\gamma+\epsilon}_{x}((1,2)\times B_1)}\leq C,\]
and thus the lemma is proved.
\qed

\noindent {\bf Proof of Proposition \ref{prop.regularityB}.~}
Exactly as in the proof of Proposition \ref{prop.regularityB}, using Lemma \ref{lem.regularity2B} in place of Lemma \ref{lem.regularity2}, we find
\[\|w\|_{C^{1+\epsilon}_{x}((1/2,1)\times B_{1/2})}\leq C.\]
Then the higher regularity of $w$ follows by considering the equation for the derivatives $D_xw$, exactly as in \cite[Theorems 6.1 and 6.2]{VDPQR}.
\qed

We can now prove Theorem \ref{Thm5}.

\noindent {\bf Proof of Theorem \ref{Thm5}.~}
As in the proof of Theorem \ref{Thm4}, by Proposition \ref{prop.regularityB} we have
\[[u]_{C^\beta_x\left(\left(t_0+\frac12r^{s(1+\frac1m)},t_0+r^{s(1+\frac1m)}\right)\times B_{r/2}(x_0)\right)}\leq Cr^{s/m-\beta}.\]
Hence, setting $\beta=k$ and $r=d(x_0)=\textrm{dist}(x_0,\partial\Omega)$, we find that
\[\bigl|D^k_xu(t,x_0)\bigr|\leq C[d(x_0)]^{\frac{s}{m}-k}\]
for all $x_0\in \Omega$ and all $t\geq t_0$.
\qed

\begin{rem}\label{schauder}
Once we know that $u\in C^\alpha$, an alternative way to show the higher regularity of $u$ (without using the results of {\rm\cite{VDPQR}}) would be the following:
If $u$ is $C^\alpha$ in a ball $(0,1)\times B_4$, then we can define $v=u\rho$ as in Lemma \ref{lem.regularity1},
and note that $v$ solves $v_t+(-\Delta)^s v^m=g\in C^\infty(B_2)$.
Then, with the same construction as in \eqref{from-w^m-to-w}-\eqref{a-alternative}, we get that $v$ solves
$v_t+Lv=g+f$ in $(0,1)\times B_1$, with $f$ and $g$ smooth in $B_{1/2}$.
Now, the $C^\alpha$ regularity of $u$ and \eqref{a-alternative} imply that the coefficients $a(t,x,y)$ are $C^\alpha$ in $x$ and $y$.
Thus, a parabolic Schauder-type estimate (like the one in \cite{CKriv}) should yield that $v$ is in fact $C^{2s+\alpha}$ in a smaller ball.
Iterating this procedure (now with $\alpha'=2s+\alpha$) one would obtain the $C^\infty$ regularity of $u$ in the $x$-variable.
We note that such a Schauder estimate has been recently established in {\rm\cite{CKriv}} for $s>\frac12$ (the results in {\rm\cite{CKriv}} hold for fully nonlinear equations and thus only for $\alpha\in(0,\bar\alpha)$, with $\bar\alpha>0$ small); see also {\rm\cite{JX}} for the case of kernels $a(x,y,t)$ which are $C^2$ in $y$, and {\rm\cite{BFV}} for the elliptic case.
Although we believe that the estimates of {\rm\cite{CKriv}} could be extended to linear equations with $s\leq\frac12$ and to all $\alpha>0$, we have preferred to use the results of {\rm\cite{VDPQR}}.
\end{rem}

Finally, we show Proposition \ref{Thm6}.

\noindent {\bf Proof of Proposition \ref{Thm6}.~}
Let $x\in K$, and let $t,s\in [t_0,T]$.
Using the equation $u_t=-(-\Delta)^s u^m$, we have the following:
\[\begin{split}
u_t(t,x)&=c_{n,s}\int_{\R^N}\frac{u^m(t,x+y)-u^m(t,x)}{|y|^{n+2s}}\,dy\\
&= c_{n,s}\int_{B_r}\frac{u^m(t,x+y)-u^m(t,x)}{|y|^{n+2s}}\,dy+c_{n,s}\int_{\R^N\setminus B_r}\frac{u^m(t,x+y)-u^m(t,x)}{|y|^{n+2s}}\,dy,
\end{split}\]
where $r\in(0,r_0)$ is a small number to be chosen later, and $r_0>0$ is a fixed number satisfying $K':=K+B_{r_0}\subset\subset \Omega$.

By Theorem \ref{Thm5}, we know that $\|u\|_{C^2_x(t_0,T]\times K')}\leq C$.
Hence, since $m>1$,
\[\left|c_{n,s}\int_{B_r}\frac{u^m(t,x+y)-u^m(t,x)}{|y|^{n+2s}}\,dy\right|\leq C\int_{B_r}\frac{|y|^2}{|y|^{n+2s}}\,dy\leq Cr^{2-2s}.\]
Therefore,
\[|u_t(t,x)-u_t(s,x)|\leq Cr^{2-2s}+C\int_{\R^N\setminus B_r}\frac{u^m(t,x+y)-u^m(s,x+y)-u^m(t,x)+u^m(s,x)}{|y|^{n+2s}}\,dy.\]
Now, by Theorem \ref{Thm4} we have that $u\in C^\gamma_t(\R^N\times [t_0,T])$, with $\gamma=\frac{1}{2m}$, which yields
\[|u_t(t,x)-u_t(s,x)|\leq Cr^{2-2s}+C\int_{\R^N\setminus B_r}\frac{|s-t|^\gamma}{|y|^{n+2s}}\,dy\leq Cr^{2-2s}+|s-t|^\gamma r^{-2s}.\]
Finally, setting $r=\min\{|s-t|^{\gamma/2},\,r_0\}$, we find
\[|u_t(t,x)-u_t(s,x)|\leq C|s-t|^\alpha,\]
with $\alpha=\gamma(1-s)$, and the Proposition is proved.
\qed

\begin{rem}
\label{rmk:alpha}
A slightly improved version of the previous argument allows one to obtain the same result with $\alpha=\min\left\{\frac{1}{2m},1-s\right\}$.
For this, one has to use in addition the interior Lipschitz regularity of $u$ in $t$ (which follows immediately from the smoothness in $x$), and split the integrals of the previous proof into $\int_{B_r}$, $\int_{B_R\setminus B_r}$, and $\int_{\R^N\setminus B_R}$, with $R>r$ chosen appropriately. We leave the details to the interested reader.
\end{rem}

\section{The case of more general domains, operators and nonlinearities}\label{sec.generalizations}

In this Section we discuss briefly how our method can be extended to more general settings. In this paper, we have preferred to stick to the prototype equation both for ease of exposition and to focus on the main ideas.

\subsection{The case of unbounded domains. }

This case includes the case of the whole Euclidean space $\RR^N$.
The extension to this case is possible by noticing that solutions $u_R(t,x)$ to the Cauchy-Dirichlet problem on a ball $B_R(x_0)\subset \Omega$ (any smaller bounded $C^{1,1}$ subdomain $B\subset\Omega$ would work) are indeed sub-solutions to the problem on the bigger domain $\Omega$, hence the instantaneous positivity result of Theorem \ref{Thm1} can be extended to the bigger domains through a simple comparison argument.
In particular, the \emph{infinite speed of propagation is valid in general domains} $\Omega\subseteq \RR^N$ (bounded or unbounded) and with no regularity assumption on $\partial\Omega$.
When the (possibly unbounded) domain $\Omega$ has a $C^{1,1}$ boundary, more precise estimates can be obtained in terms of the distance to the boundary, as in Theorem \ref{Thm2}.

As a consequence of the comments above, the $C^\infty_x$ interior regularity hold true for all domains.
When the domain has $C^{1,1}$ boundary, then the $C^{s/m}_x$ regularity up to the boundary follows exactly as in Theorem \ref{Thm4}.
Notice that the results in Propositions \ref{prop.regularity} and \ref{prop.regularityB} are local estimates, in the sense that only require the equation to be satisfied in a ball, and the function $u$ to be bounded in $\RR^N$.
This means that both results can be applied directly to any bounded solution $u$ (actually, using smoothing effects,
one can consider even more singular initial data), even when the domain $\Omega$ is not bounded or not regular.

\subsection{More general kernels}

We can consider a family of nonlocal operators with more general kernels than the one of the fractional Laplacian, as suggested in formula \eqref{eq:general L}\,,
\begin{equation}\label{eq:general L.2}
\A f(x):={\rm PV}\int_{\RR^N}\big(f(x)-f(y)\big)\,\K(x,y)\dy\,.
\end{equation}
The typical assumption on the kernels is $\K(x,y)=\K(y,x)$ together with
\begin{equation}\label{kernels}
\frac{\lambda}{|x-y|^{n+2s}}\leq \K(x,y)\leq \frac{\Lambda}{|x-y|^{n+2s}},
\end{equation}
being $\lambda,\Lambda$ positive constants.
In literature this type of operators are sometimes called ``rough kernels''.

In order our positivity result to be true, we do not need any extra assumption on the kernel, and the proof works for any kernel of the form \eqref{kernels}.
When $u$ and the kernels \eqref{kernels} are regular enough, then the proof is essentially the same.
In order to justify rigorously the argument for a general operator \eqref{eq:general L.2}, one should use the weak formulation of the equation.
If we want the sharp two-sided estimate from Theorem~\ref{Thm2}, then for the upper bound one needs to ensure the validity of suitable Green function estimates on bounded domains,
namely
\begin{equation}\label{Green.est.general.K}
G(x,y)\asymp \frac{1}{|x-y|^{N-2s}}
\left(\frac{d(x)^s}{|x-y|^s}\wedge 1\right)
\left(\frac{d(y)^s}{|x-y|^s}\wedge 1\right)\,.
\end{equation}
This estimate is known to hold for the fractional Laplacian in bounded $C^{1,\alpha}$ domains \cite{Kim-Coeff}, but do not hold for general kernels \eqref{kernels}.
Still, for general operators with kernels \eqref{kernels}, one can obtain an $L^\infty$ bound for the solution, which combined with our lower bound would imply that the solution is bounded between two positive constants in any compact subset $K\subset\subset \Omega$.
This, combined with our regularity arguments, would yield the $C^{0,\alpha}$ interior regularity of solutions.

\noindent\textbf{Results for more general kernels. }We summarize the results for the case of more general kernels, in the case when $\Omega$ is a bounded domain. The case of $\Omega$ unbounded is also possible as explained in the previous section.
\begin{itemize}[leftmargin=*]\itemsep2pt \parskip3pt \parsep0pt
\item \textit{Positivity and infinite speed of propagation. } For general kernels as in \eqref{kernels}, we have \textit{infinite speed of propagation } for all $0<s<1$.
    The sharp two sided bound from Theorem \ref{Thm2} holds whenever the Green function estimate \eqref{Green.est.general.K} is true.

\item \textit{H\"older regularity. } For general kernels as in \eqref{kernels} we have an interior space-time H\"older regularity estimate, but not up to the boundary.
    When \eqref{Green.est.general.K} holds, then we obtain the \textit{sharp $C^{s/m}_x$ regularity at the boundary}, exactly as in Theorem \ref{Thm4}.

\item \textit{Higher regularity. }To prove the higher regularity of solutions we have used the results of \cite{VDPQR}, which are true only for $(-\Delta)^s$.
    Still, when the kernels are regular in $x$ and $y$, our methods work
    provided one can use parabolic Schauder estimates; see Remark \ref{schauder}.
    In case $s>1/2$, these Schauder estimates follow from the results of \cite{CKriv}.

\end{itemize}

\begin{rem}\label{rem.PQRV}
After the writing of this paper was completed, we learned that,
at the very same time and independently of us,
de Pablo, Quiros, and Rodriguez proved in {\rm\cite{DPQR}} the H\"older regularity for solutions to the Cauchy problem in the whole Euclidean space $\RR^N$ for the the equation $u_t=-\A\n(u)$\,, when $\A$ is an operator with rough kernels (as the first example above) and $\n$ is a general nonlinearity.
Their result is based on a nonlocal version of the De Giorgi method,
and generalizes in  many   aspects the previous regularity results of {\rm\cite{AC2010,DPQRV2,VDPQR}}.
As explained above, our method could be adapted to cover the case of rough kernels on the whole space, as in {\rm\cite{DPQR}}\,, and in the case of the Cauchy-Dirichlet problem on domains.
\end{rem}

\subsection{More general nonlinearities}
We can allow more general nonlinearities than the pure power case $\n(u)=|u|^{m-1}u$ with $m>1$.
Indeed, following the general setup given in \cite{BV-PPR2-1}\,, we can allow for continuous and non-decreasing functions $\n:\RR\to\RR$, with the  normalization $F(0)=0$, satisfying the condition
\begin{enumerate}
\item[(N1)] $\n\in C^1(\RR\setminus\{0\})$, $\n/\n'\in {\rm Lip}(\RR)$, and there exist $\mu_0,\mu_1\in (0,1)$ such that
\[
1-\mu_1\le \left(\frac{\n}{\n'}\right)'\le 1-\mu_0\,,
\]
where $\n/\n'$ is understood to vanish if $\n(r)=\n'(r)=0$ or $r=0$\,.
\end{enumerate}
The main example is $\n(u)=|u|^{m-1}u$, with $m>1$ (in which case $\mu_1=\mu_0=(m-1)/m$)\,, and a simple variant is the combination of two powers, so that one of them gives the behavior near $u=0$, the other one the behavior near $u=\infty$.
A quite complete theory of existence, uniqueness, and a priori estimates, including the GHP for large times, for weak dual solutions to the Cauchy-Dirichlet problem for the equation $u_t=-\A \n(u)$ under certain assumptions on the kernels has been given in \cite{BV-PPR2-1, BV-PPR2-2}, see also Appendix \ref{App.I} for some more details.
The key point of assumption $(N1)$ is that it guarantees monotonicity in time estimates, namely the function $t\mapsto t^{1/\mu_0}\n(u(t,x))$ is non-decreasing for almost all $x\in\Omega$\,, as it has been proven by Crandall and Pierre in \cite{CP-JFA}.

\noindent\textbf{Results for more general nonlinearities. }We summarize the results for the case of more general nonlinearities, in the case when $\Omega$ is a bounded domain. The case of $\Omega$ unbounded could also treated as explained above.
\begin{itemize}[leftmargin=*]\itemsep2pt \parskip3pt \parsep0pt
\item \textit{Positivity and infinite speed of propagation. }A minor modification of proof of Theorem \ref{Thm1} allows one to prove the positivity result for small times also in this case; we can also obtain Global and local Harnack inequalities in the spirit of Theorems \ref{Thm2} and \ref{Thm3} (in a slighlty modified form) by combining the results of this paper with the ones of \cite{BV-PPR2-2}.
    In particular, we still have \textit{infinite speed of propagation.}

\item \textit{$C^\alpha_{x,t}$ regularity. } As a consequence, the $C^\alpha$ regularity results hold both in bounded and unbounded domains and also for operators with rough kernels, generalizing the recent results of \cite{DPQR} and the previous results of \cite{AC2010} and \cite{VDPQR}.

\item\textit{Higher regularity. }Finally, when $\A=(-\Delta)^s$, the higher regularity results of Theorems \ref{Thm5} and \ref{Thm6} hold for general nonlinearities $\n$ under the running assumptions both in bounded and unbounded domains; when $\A$ is a Levy operator with rough kernels, higher regularity shall follow by Schauder estimates, as already commented in the previous section.
\end{itemize}

\subsection{The classical case $s=1$}

The classical case has been intensively studied since the 1980's, by many authors \cite{Ar69,Ar70b,ArCaVa85,CaFr79b, CaFr80,CaVaWo87,CaWo90,DaHaLe01,DaHa98,DiB86,DiB88,DGVacta,GiPe81}, see also the books \cite{DaskaBook, DiBook,DGVbook,VazBook}, as already mentioned in the introduction.
We refer to \cite[Chapter 19]{VazBook} for an exposition about the regularity results for the PME.

We first recall that the positivity result for small times is not valid for $s=1$, since when $s=1$ there is finite speed of propagation.
In any case, after the waiting time $t_*$, the GHP of Theorem \ref{Thm2} holds true, cf. \cite{BV-PPR1, BV-PPR2-2} and as a consequence the regularity results of Theorems \ref{Thm4}, \ref{Thm5} and \ref{Thm6} hold also in the limit case $s=1$.

Since the equation is local, it follows by parabolic regularity that, for $t \geq t_*$, the solution is $C^\infty$ in space-time
in the interior of $\Omega$. However,
to our best knowledge, even in this local case higher order boundary regularity is not known.

\section{Appendix I. Weak dual solutions and their basic properties}\label{App.I}

The fractional Laplacian operator that we are considering in this paper is defined through the singular integral representation in the whole space
\begin{equation}\label{sLapl.Rd.Kernel}
(-\Delta)^{s}  g(x)= c_{N,s}\mbox{
P.V.}\int_{\mathbb{R}^N} \frac{g(x)-g(z)}{|x-z|^{N+2s}}\,dz,
\end{equation}
where $c_{N,s}>0$ is a normalization constant.
Since we study the homogeneous Dirichlet problem, we ``restrict'' the operator to functions that are zero outside a bounded domain $\Omega$
of class $C^{1,1}$.
It is worth mentioning that in \cite{BV-PPR1,BV-PPR2-1,BV-PPR2-2,BSV2013} this operator has been called \textit{restricted fractional Laplacian}, RFL for short, in order to distinguish it from other possible (non-equivalent) definitions of the Dirichlet fractional Laplacian on domains.
In this case, the initial and boundary conditions associated to the fractional diffusion equation $u_t+(-\Delta)^s u^m=0$ are
\begin{equation}\label{Bound.Cond.RFPME}
\left\{
\begin{array}{lll}
u(t,x)=0\; &\mbox{in }(0,\infty)\times\RR^N\setminus \Omega\,,\\
u(0,\cdot)=u_0\; &\mbox{in }\Omega\,.
\end{array}
\right.
\end{equation}
We refer to \cite{BSV2013} for a careful construction of the RFL in the framework of fractional Sobolev spaces.
The operators $(-\Delta)^s$, $s\in(0,1)$, are infinitesimal generators of stable and radially symmetric L\'evy processes \cite{Levy,Levy2}.
These stochastic processes play an important role in Probability, and have been used in the last years to model prices in Finance \cite{finance,finance2}, or anomalous diffusions in Physics \cite{phys2,phys}, Biology and Ecology \cite{nature,nature2}, among others.

So defined, $\A=(-\Delta)^s$ is a self-adjoint operator on $\LL^2(\Omega)$\,, with a discrete spectrum: we denote by $\lambda_{s, j}>0$, $j=1,2,\ldots$ its eigenvalues written in increasing order and repeated according to their multiplicity, and denote by $\{\phi_{s, j}\}_j$ the corresponding set of eigenfunctions, normalized in $L^2(\Omega)$.
The corresponding eigenfunctions are known to be only H\"older continuous up to the boundary cf. \cite{RS-Dir, Grub1}.
We will denote by $\Phi_1$ the first positive eigenfunction, and we recall that
\begin{equation}\label{Phi.Estimates}
\Phi_1(x)\asymp \dist(x, \R^N\setminus\Omega)^s;
\end{equation}
see \cite{RS-Dir}.

The inverse of $\A=(-\Delta)^s$ in $\Omega$, $\AI:\LL^1(\Omega)\to\LL^1(\Omega)$, can be defined by means of the Green function of $\A$ as follows
\[
\AI f(x):=\int_\Omega G(x,y)f(y)\dy
\]
Moreover, we know by the results of \cite{Chen-Song,Kul} that, when $N>2s$,
\begin{equation}\label{Green.Estimates}
G(x,y)\asymp \frac{1}{|x-y|^{N-2s}}
\left(\frac{\dist(x, \partial\Omega)}{|x-y|}\wedge 1\right)^s
\left(\frac{\dist(y, \partial\Omega)}{|x-y|}\wedge 1\right)^s\,.
\end{equation}
The above facts allow one to use the theory of existence and uniqueness and the estimates developed in \cite{BV-PPR1,BV-PPR2-1,BV-PPR2-2,BSV2013}.
For the convenience of the reader, we now briefly rescall them, referring to those papers for further details and generalizations.

As shown in \cite{BV-PPR1,BV-PPR2-1,BV-PPR2-2,BSV2013}, the results in this section hold both for more general nonlinearities $F$ other than $|u|^{m-1}u$\,, and for more general operators $\A$.  However, in order to avoid unnecessary confusion, we restrict ourselves to the case of nonnegative solutions and to the power case, hence here $F(u)=|u|^{m-1}u=u^m$.
We finally recall that $s\in (0,1]$, $N>2s$, and $m>1$.

\begin{defn}\label{Def.Very.Weak.Sol.Dual} A function $u$ is called a {\sl weak dual} solution to the Dirichlet Problem for the equation $u_t=-\A u^m$ in $Q_T=(0,T)\times \Omega$ if:
\begin{itemize}
\item $u\in C((0,T): \LL^1_{\Phi_1}(\Omega))$\,, $u^m \in \LL^1\left((0,T):\LL^1_{\Phi_1}(\Omega)\right)$;
\item  The identity
\begin{equation}
\displaystyle \int_0^T\int_{\Omega}\AI (u) \,\dfrac{\partial \psi}{\partial t}\,\dx\dt
-\int_0^T\int_{\Omega} u^m\,\psi\,\dx \dt=0.
\end{equation}
holds for every test function $\psi$ such that  $\psi/\Phi_1\in C^1_c((0,T): \LL^\infty(\Omega))$\,.
\end{itemize}\end{defn}
We are interested in solving the Cauchy-Dirichlet problem, consisting of the equation $u_t=-\A u^m$ with homogeneous Dirichlet conditions plus given initial data.
\begin{defn}\label{Def.Very.Weak.Sol.Dual.CDP}
A {\sl weak dual} solution to the Cauchy-Dirichlet problem is a weak dual solution to the Dirichlet problem for the equation $u_t=-\A u^m$ satisfying $u\in C([0,T): \LL^1_{\Phi_1}(\Omega))$  and $u(0,x)=u_0\in \LL^1_{\Phi_1}(\Omega)$.
\end{defn}
This kind of solution has been first introduced in \cite{BV-PPR1,BV-PPR2-1}. Roughly speaking, we are considering the weak solution to the ``dual equation'' $\partial_t U=- u^m$\,, where $U=\AI u$\,, posed on the bounded domain $\Omega$ with homogeneous Dirichlet conditions. The advantage of considering the dual problem is that the boundary conditions are implicitly contained in $\AI$.

\begin{thm}[Existence and Uniqueness of weak dual solutions \cite{BV-PPR1,BV-PPR2-1}]\label{thm.exist.uniq.weak.dual}
For every  nonnegative $ u_0\in\LL^1_{\Phi_1}(\Omega)$ there exists a unique minimal weak dual solution to the Cauchy-Dirichlet problem for the equation $u_t=-\A u^m$. Such a solution is obtained as the monotone limit of the semigroup (mild) solutions that exist and are unique. We call it the minimal solution. The minimal weak dual solution is continuous in the weighted space $u\in C([0,\infty):\LL^1_{\Phi_1}(\Omega))$, and the standard comparison result holds.
\end{thm}
\noindent\textbf{Remarks. }(i) The construction of the minimal weak dual solution is realized by approximation of the initial datum $u_0\in \LL^1_{\Phi_1}(\Omega)$ from below, in terms of
a sequence $u_{0,n}\in\LL^1$. This sequence generates a sequence of unique mild (semigroup) solutions as it is shown in \cite{CP-JFA}. We refer to \cite{BV-PPR2-1} for further details and for the proof of the fact that mild solutions are indeed weak dual solutions.

\noindent(ii) A large class of solutions fall into this class of solutions: weak, mild, strong, $H^*$, see \cite{BV-PPR1, BV-PPR2-1, BSV2013}. Weak dual solutions satisfy a quite complete set of a priori estimates, namely absolute upper bounds (sharp up to the boundary), smoothing effects (instantaneous and even backward in time), positivity estimates (with sharp boundary behaviour). Such bounds combine to yield the following sharp global Harnack-type estimates.
\begin{thm}[Global Harnack Principle \cite{BV-PPR1, BV-PPR2-2}]\label{thm.GHP.PME}
There exist constants $H_0, H_1, c_*>0$ such that setting,
\begin{equation}\label{thm.GHP.PME.0}
t_*= \frac{c_*}{\left(\int_{\Omega}u_0\Phi_1\dx\right)^{m-1}}\,,
\end{equation}
the following inequality holds for all $t\ge t_*$ and all $x\in \Omega$:
\begin{equation}\label{thm.GHP.PME.1}
H_0\,\frac{\Phi_1(x)^{\frac{1}{m}}}{t^{\frac{1}{m-1}}} \le \,u(t,x)\le H_1\, \frac{\Phi_1(x)^{\frac{1}{m}}}{t^{\frac{1}{m-1}}}
\end{equation}
\end{thm}
We recall that  the upper bound in formula \eqref{thm.GHP.PME.1} holds true for all times $t>0$, while the lower bound only holds for $t\ge t_*$. The first main result of the present paper is devoted to prove an analogous lower bound for small times; this problem was left open in \cite{BV-PPR1,BV-PPR2-2}.

As a corollary of the above estimates we also have local forms of Harnack inequalitites:
\begin{cor}[Local Harnack Inequalities of Backward Type \cite{BV-PPR1, BV-PPR2-2}]\label{thm.Harnack.Local.backward}
Under the assumptions of Theorem \ref{thm.GHP.PME.0}, there exist constants $H_2$, $L_0>0$ such that the following inequality holds
for all $t\ge t_*$ and all $B_R(x_0)\in \Omega$:
\begin{equation}\label{thm.Harnack.PME.2}
\sup_{x\in B_R(x_0)}u(t,x)\le H_2\, \inf_{x\in B_R(x_0)}u(t+h,x)\qquad\mbox{for all $0\le h\le t_*$\,.}
\end{equation}
\end{cor}
\noindent\textbf{Remark. }All the above constants $H_2,H_1,H_0, c_*>0$ depend only on $N, m, s,$ and $\Omega$\,,  and have an explicit form given in \cite{BV-PPR1, BV-PPR2-1,BV-PPR2-2}.

\section{Appendix II. Approximate solutions}\label{appendix.approximate}

Let us fix $\delta>0$.
To prove our lower bounds in Section \ref{sect:lower},
we considered the ``larger'' approximate problem:
\begin{equation}\label{Dirichlet.Udelta}
\left\{
\begin{array}{lll}
\partial_t u_\delta=-\A u_\delta^m & \qquad\mbox{for any $(t,x)\in (0,\infty)\times \Omega$}\\
u_\delta(t,x)=\delta &\qquad\mbox{for any $(t,x)\in (0,\infty)\times (\RR^N\setminus\Omega)$}\\
u_\delta(0,x)=u_0(x)+\delta &\qquad\mbox{for any $x\in\Omega$}\,.
\end{array}
\right.
\end{equation}
Letting now $u_\delta=v_\delta+\delta$, we see that $v_\delta$ solves
\begin{equation}\label{Dirichlet.Vdelta}
\left\{
\begin{array}{lll}
\partial_t v_\delta=-\A\left[(v_\delta+\delta)^m -\delta^m \right]& \qquad\mbox{for any $(t,x)\in (0,\infty)\times \Omega$}\\
v_\delta(t,x)=0 &\qquad\mbox{for any $(t,x)\in (0,\infty)\times (\RR^N\setminus\Omega)$}\\
v_\delta(0,x)=u_0(x) &\qquad\mbox{for any $x\in\Omega$}\,.
\end{array}
\right.
\end{equation}
which is a Dirichlet problem with a nonlinearity $\n_\delta(v)=(v+\delta)^m -\delta^m$\,, which is neither singular nor degenerate at $v=0$, since $F_\delta(0)=0$ and $F_\delta'(0)=m\delta^{m-1}>0$. Then, the theory developed in \cite{BSV2013, BV-PPR1, BV-PPR2-1, BV-PPR2-2} allows one to show that $v_\delta$ enjoys a number of properties and estimates that we collect in this Appendix for reader's convenience.

\subsection{Existence and Uniqueness}
The definition of solution that we use for the approximate Dirichlet problem \eqref{Dirichlet.Vdelta} is the same as in Definition \ref{Def.Very.Weak.Sol.Dual.CDP}. The approximate solutions $v_\delta$  will always be the minimal weak dual solution of the Cauchy-Dirichlet problem \eqref{Dirichlet.Vdelta}. Let us state an existence result, which follows by the results of \cite{BSV2013} applied to the present case.
The fractional Sobolev spaces will be denoted by $H^*$, where $H^*=H^{-s}(\Omega)=(H^{s}_{0})^*$ for any $s\ne 1/2$ and $H^*=(H^{1/2}_{00})^*$ when $s=1/2$.

\begin{thm}[Existence and uniqueness of approximate solutions]\label{exist.uniq.approx.sols}
For every $0\le u_0\in \LL^1_{\Phi_1}(\Omega)$ there exists a unique global weak dual solution $v_\delta$ of  Problem \eqref{Dirichlet.Vdelta}.
In addition the solution map $S_t: u_0\mapsto v_\delta(t)$ defines a semigroup of (non-strict) contractions in $H^*(\Omega)$, i.\,e.,
\begin{equation}\label{contractivity.H*}
\|v_\delta(t)-\bar{v}_\delta(t)\|_{H^*(\Omega)}\le\|v_\delta(0)-\bar{v}_\delta(0)\|_{H^*(\Omega)},
\end{equation}
which turns out to be also compact in $H^*(\Omega)$.
\end{thm}
\noindent\textbf{Remark. }An analogous existence and uniqueness result can be stated for $u_\delta=v_\delta+\delta$.

\noindent {\sl Proof.~}As mentioned above, existence and uniqueness of weak dual solutions for the approximate problem \eqref{Dirichlet.Vdelta} follows by the results of \cite{BSV2013}. More precisely, in \cite[Theorem 2.2]{BSV2013} the authors prove existence and uniqueness for a different class of solutions from the one we consider here: the so-called $H^*$-solutions. Porous media-type equations are shown to generate a nonlinear semigroup of contraction in $H^*$, and existence and uniqueness of strong $H^*$-solutions are established\,, for any $u_0\in H^*$. This is not enough for our purpose here, as we want to prove existence and uniqueness of weak dual solutions with $0\le u_0\in \LL^1_{\Phi_1}$, which is a different class of solutions.

\noindent Recall that nonnegative functions of $H^*$ belong to $\LL^1_{\Phi_1}$. Indeed, by Cauchy-Schwartz inequality,
\begin{equation}\label{L1.H}
\|f\|_{\LL^1_{\Phi_1}(\Omega)}=\int_\Omega f\Phi_1\dx=\int_\Omega \AIM f \AM\Phi_1\dx\le \|f\|_{H^*}\|\Phi_1\|_{H}\le \lambda_1\|f\|_{H^*}\,.
\end{equation}
Viceversa, if $0\le f\in \LL^1_{\Phi_1}$ and $f/\Phi_1\in \LL^\infty$\,, then $f\in H^*$:
\begin{equation}\label{L1.H.2}
\|f\|_{H^*}= \int_\Omega f\AI f\dx \le \left\|\frac{f}{\Phi_1}\right\|_{\LL^\infty(\Omega)} \int_\Omega \Phi_1\AI u
= \lambda_1^{-1} \left\|\frac{f}{\Phi_1}\right\|_{\LL^\infty(\Omega)}\|f\|_{\LL^1_{\Phi_1}(\Omega)}\,.
\end{equation}
Notice that if $f$ is compactly supported in $\Omega$\,, we always have $f/\Phi_1\in \LL^\infty$, since $\Phi_1\asymp \dist(\cdot,\partial\Omega)^s$.

Thanks to these observation, we can
approximate from below $0\le u_0\in \LL^1_{\Phi_1}$ by a sequence of nonnegative compactly supported functions $0\le u_{0,n}\in C_c^\infty\subset H^*\cap \LL^1_{\Phi_1}$, with $u_{0,n}\to u_0$ in the strong $\LL^1_{\Phi_1}$ topology. Let $u_n(t)$ be the corresponding unique $H^*$ solutions, also belonging to $H^*\cap\LL^1_{\Phi_1}$ by \eqref{L1.H}.
Next, we recall \cite[Lemma 5.3]{BSV2013} which shows that nonnegative $H^*$-solutions are weak dual solutions. The comparison principle holds for $H^*$ solutions, hence we have that, for any fixed $t>0$, the monotone sequence $0\le u_{n}(t)\le u_{n+1}(t)$ converges pointwise to $u(t)\ge 0$, defined as the limit of the monotone sequence $u_n(t,x)$ as $n\to \infty$. The solution $u$ constructed in this way is called minimal weak dual solution in \cite{BV-PPR1,BV-PPR2-1} (the only difference in the construction is the class of approximating solutions, here we use $H^*$ instead of mild solutions). A standard limiting process shows that also $u$ is a weak dual solution, the key point being that $u-u_n\ge 0$, which implies
\[
\|u(t)-u_n(t)\|_{\LL^1_{\Phi_1}(\Omega)}\le \|u_0-u_{0,n}\|_{\LL^1_{\Phi_1}(\Omega)}\,
\]
cf. Lemma \ref{Lem.0} for a proof. Finally, uniqueness of the minimal weak dual solution $u$ can be proven exactly as in \cite[Theorem 4.5]{BV-PPR2-1}.\qed

\subsection{Comparison results involving approximate solutions}

We prove now some comparison results implying
\eqref{approx.soln.positivity}, \eqref{approx.soln.comparison1}, and \eqref{approx.soln.comparison}.
Recall that these properties are crucial in the proof of the lower bounds of Theorem \ref{Thm1}. We will prove them in the two following lemmata.
\begin{lem}
Let $\delta > 0$, let $v_\delta$ be a solution of \eqref{Dirichlet.Vdelta}, and let $u$ be a weak dual solution, both corresponding to the same initial datum $u_0\ge 0$. Then \eqref{approx.soln.comparison} holds.
\end{lem}

\noindent {\sl Proof.~}
Notice that the statement is equivalent to showing that $v_\delta(t,x)+\delta\ge u(t,x)$.

The formal proof is as follows: let $H$ be the Heaviside function, and set $(u)_+:=\max\{u,0\}$.
Then, since
\[
\A\left[u^m-(v_\delta+\delta)^m+\delta^m\right]=\A\left[u^m-(v_\delta+\delta)^m\right]
\]
and
\[
H(u-\delta-v_\delta)=H\left[u^m-(v_\delta+\delta)^m\right],
\]
for all $0\leq t_0 \leq t$ we have
\begin{equation}\label{formal.comparison}\begin{split}
\int_{t_0}^t\int_\Omega \partial_t(u-\delta-v_\delta)_+ \dx\rd\tau
&=-\int_{t_0}^t\int_\Omega H(u-\delta-v_\delta) \A\left[u^m-(v_\delta+\delta)^m+\delta^m\right]\dx\rd\tau\\
&=-\int_{t_0}^t\int_\Omega H\left[u^m-(v_\delta+\delta)^m\right] \A\left[u^m-(v_\delta+\delta)^m\right]\dx\rd\tau\le 0,
\end{split}
\end{equation}
where last inequality follows by the accretivity of $\A$ and the fact that $H$ is monotone nondecreasing.
More precisely, recall here the Stroock-Varopoulos inequality,
stating that for all smooth increasing functions  $\phi:\RR\to\RR$ one has
\begin{equation}\label{Stroock.Varopoulos.full}
\int_\Omega \phi(f(x)) \A f(x) \dx\ge \int_\Omega \left| \AM \eta(f(x))\right|^2 \dx\ge 0,
\end{equation}
where $\eta:\RR\to\RR$ is such that $(\eta')^2=\phi'$.
By a simple approximation argument, we can take $\phi=H$, which proves the last inequality in
\eqref{formal.comparison}.

For the convenience of the reader, we recall that
the Stroock-Varopoulos follows by the following simple argument: first of all we observe that,  for all $a,b\ge 0$,
\[
(a-b)(\phi(a)-\phi(b))\ge (\eta(a)-\eta(b))^2\,,\qquad\mbox{where $(\eta')^2=\phi'$}\,.
\]
Then, thanks to this inequality and using the fact that $\A$ can be represented by a symmetric nonnegative kernel $\K(x,y)$, we get
\[\begin{split}
\int_\Omega \phi(f(x)) \A f(x) \dx&=\int_\Omega \int_{\RR^N}\phi(f(x))\,[f(x)-f(y)]\K(x,y)\dy\dx\\
&=\frac{1}{2}\int_\Omega \int_{\RR^N}[\phi(f(x))-\phi(f(y))]\,[f(x)-f(y)]\,\K(x,y)\dy\dx\\
&\ge \frac{1}{2}\int_\Omega \int_{\RR^N}[\eta(f(x))-\eta(f(y))]^2\,\K(x,y)\dy\dx\ge 0.
\end{split}
\]
One then obtains the Stroock-Varopoulos inequality \eqref{Stroock.Varopoulos.full} once one observes that the last term of the above inequality can also be rewritten as
\[
\frac{1}{2}\int_\Omega \int_{\RR^N}[\eta(f(x))-\eta(f(y))]^2\,\K(x,y)\dy\dx= \int_{\RR^N} \eta(f)\A \eta(f)\dx
=\int_{\RR^N}\left|\AM \eta(f)\right|^2\dx
 \ge\int_\Omega\left|\AM \eta(f)\right|^2\dx\,.
\]
This argument shows that the $L^1$ norm of $(u-\delta-v_\delta)_+$ is decreasing in time, hence it must be identically zero as it vanishes for $t=0$.

Although this proof is formal, it can be made rigorous by taking a smooth approximation of the function $(\cdot)_+$ and working with
weak energy solutions (instead of weak dual solutions). We refer to \cite{BV-PPR1} for more details.
Also, we mention that another proof can be obtained by using the extension problem for the fractional Laplacian, as it has been done in \cite{DPQRV2}.\qed

This same argument can be used to show also the validity of \eqref{approx.soln.comparison1}.
We finally prove a general comparison principle that, as an immediately consequence, yields \eqref{approx.soln.positivity}.

\begin{lem}[Comparison for approximate solutions]
Let $\delta> 0$ and let $v_\delta,w_\delta$ be two solutions of \eqref{Dirichlet.Vdelta}. Then for all $t\ge t_0\ge 0$ we have
\begin{equation}\label{comparison.approx.soln}
\int_\Omega(v_\delta(t,x)-w_\delta(t,x))_+\,\dx\le \int_\Omega(v_\delta(t_0,x)-w_\delta(t_0,x))_+\,\dx\,.
\end{equation}
As a consequence, if $w_\delta(0)\le v_\delta(0)$ then $w_\delta\le v_\delta$ in $(0,\infty)\times\Omega$.
\end{lem}

Thanks to this result, if $u_0\ge 0$ then $v_\delta\ge 0$\,, that is $u_\delta\ge \delta$\,.

\noindent {\sl Proof.~}The proof is an adaptation from the proof for the case $s=1$, see e.g. \cite{VazBook}, or also \cite{DPQRV2} where the same result is proved for solutions to the Cauchy problem in the case $\delta=0$. Therefore we just sketch the proof for convenience of the reader.

Let $H$ be the Heaviside function and let $(u)_+:=\max\{u,0\}$. Then, since $H(v_\delta(t,x)-w_\delta(t,x))=H\left[(v_\delta+\delta)^m-(w_\delta+\delta)^m)\right]$ and
using \eqref{Stroock.Varopoulos.full} with $f=(v_\delta+\delta)^m-(w_\delta+\delta)^m$, we get
\[\begin{split}
\int_{t_0}^t\int_\Omega \partial_t(v_\delta(t,x) -w_\delta(t,x))_+ & \dx\rd\tau
=-\int_{t_0}^t\int_\Omega H(v_\delta(t,x)-w_\delta(t,x)) \A\left[(v_\delta+\delta)^m - (w_\delta+\delta)^m\right]\dx\rd\tau\\
&=-\int_{t_0}^t\int_\Omega H\left[(v_\delta+\delta)^m-(w_\delta+\delta)^m\right] \A\left[(v_\delta+\delta)^m - (w_\delta+\delta)^m\right]\dx\rd\tau\le 0.
\end{split}\]
As before, the proof can be made rigorous by means of smooth approximation of $H$ and $(\cdot)_+$\,, and an approximation of weak dual solutions by weak energy solutions.\qed

\subsection{Boundedness and regularity}\label{app.reg.approx}

A small modification of the proof of \cite[Theorem 2.2]{BV-PPR2-1} or \cite[Theorem 5.1]{BV-PPR1} allows us to prove absolute upper bounds for the approximate solutions $u_\delta$.
Indeed, by the Benilan-Crandall estimates \cite{BCr}, the following inequality holds in the distributional sense:
\begin{equation}\label{monotonicity.approx}
\partial_t u_\delta(t,x)\ge -\frac{u_\delta(t,x)}{(m-1)t}
\end{equation}
As a consequence, the function $t\mapsto t^{\frac{1}{m-1}}u_\delta(t,\cdot)$ is monotone non decreasing for $t>0$. We then transfer the monotonicity of $u_\delta$ to $v_\delta=u_\delta-\delta$ to obtain the following:

\begin{prop}[Absolute upper estimate for approximate solutions]\label{thm.Upper.Approximate}
Let  $v_\delta$ be an approximate solution to the Cauchy-Dirichlet problem \eqref{Dirichlet.Vdelta}.
Then, there exists a constant $K_1>0$ such that the following estimates hold true
$$
 \|v_\delta(t)\|_{\LL^\infty(\Omega)}\le\, K_1\,t^{-\frac{1}{m-1}}+\,\delta\,,\qquad\qquad\mbox{for all } \ t> 0\,.
$$
\end{prop}
The constant $K_1$ depends only on $N, m, s$ and $\Omega$, and has an explicit form given in \cite{BV-PPR1, BV-PPR2-1}.

\noindent\textit{Regularity for $u_\delta$. }Once we know that $v_\delta$ are nonnegative and bounded, we can see that $u_\delta$ is strictly positive and bounded: more precisely, for almost all $(t,x)\in (0,\infty)\times \Omega$
we have
\[
\delta\le u_\delta(t,x)\le K_1\,t^{-\frac{1}{m-1}}+\,2\delta
\]
Therefore, we can apply to the approximate solutions $u_\delta$ the same arguments of the proofs of Theorems \ref{Thm4}, \ref{Thm5}, and \ref{Thm6}, to deduce the same regularity results for $u_\delta$; the approximate solution $u_\delta$ turns out to be globally H\"older continuous in space-time.
Moreover, they are classical solution in the interior. This completely justifies all the computations in the proofs involving approximate solutions.

\vspace{5mm}

\noindent {\textbf{\large \sc Acknowledgments.} M.B. has been partially funded by Project MTM2011-24696 and MTM2014-52240-P (Spain).
A.F. is supported by NSF Grants DMS-1262411 and
DMS-1361122.
M.B. would like to thank the Mathematics Department of Texas University at Austin for its kind hospitality, where part of this work has been done.

\vskip .3cm



\addcontentsline{toc}{section}{~~~References}
\begin{thebibliography}{00}
\small

\bibitem{Ar69}D. G. Aronson, \textit{Regularity propeties of flows through porous media. }SIAM J. Appl. Math. \textbf{17} (1969), 461-–467.

\bibitem{Ar70b}D. G. Aronson, \textit{Regularity properties of flows through porous media: A counterexample. }SIAM J. Appl. Math. \textbf{19} (1970) 299--307.

\bibitem{ArCaVa85}D.G. Aronson, L.A. Caffarelli, J.L. V\'azquez, \textit{Interfaces with a corner-point in one-dimensional porous medium flow, }Comm. Pure Appl. Math. \textbf{38} (1985), 375--404.

\bibitem{Ar-Pe} D. G. Aronson, L. A. Peletier. \textit{Large time behaviour of solutions of the porous medium equation in bounded domains}, J. Differential
Equations \textbf{39} (1981), 378--412.

\bibitem{AC2010} { I. Athanasopoulos, L. A. Caffarelli, }\emph{Continuity of the temperature in boundary heat control problems}. {Adv. Math. } \textbf{224} (2010), 293--315.


\bibitem{BFV} B. Barrios, A. Figalli, E. Valdinoci,
\emph{Bootstrap regularity for integro-differential operators, and its application to nonlocal minimal surfaces}.
Ann. Sc. Norm. Super. Pisa Cl. Sci. (5), \textbf{13} (2014), no. 3, 609--639.

\bibitem{BCr}{\rm  P. B\'enilan, M.~G. Crandall.} \textit{Regularizing effects of homogeneous evolution equations,} {\rm
Contributions to Analysis and Geometry} (suppl.  to Amer. J. Math.), Johns Hopkins Univ. Press, Baltimore, Md., 1981. 23-39.

\bibitem{Levy2} J. Bertoin, \emph{L\'evy Processes}, Cambridge University Press, Cambridge, 1996.

\bibitem{BGT2000} {\rm M. Bologna,  P. Grigolini, C. Tsallis, }\textit{Anomalous diffusion associated with nonlinear fractional derivative Fokker-Planck-like equation: Exact time-dependent solutions}, \textrm{Physical Review E }\textbf{62}, (2000).

\bibitem{BSV2013} M. Bonforte, Y. Sire, J.~L. V\'azquez, \textit{Existence, Uniqueness and Asymptotic behaviour for fractional porous medium on bounded domains.  Relation with elliptic equations}. To appear in Disc. Cont. Dyn. Syst. (2015).

\bibitem{BV-PPR1} M. Bonforte,  J. L. V\'azquez, \emph{A Priori Estimates  for Fractional Nonlinear  Degenerate Diffusion Equations on bounded domains},  Arch. Rat. Mech. Anal. \textbf{218} (2015), no. 1, 317--362. 

\bibitem{BV-PPR2-1} M. Bonforte,  J. L. V\'azquez, \emph{Fractional Nonlinear Degenerate Diffusion Equations on Bounded Domains Part I. Existence, Uniqueness and Upper Bounds. }Nonlin. Anal. TMA \textbf{131} (2016), 363--398.

\bibitem{BV-PPR2-2} M. Bonforte,  J. L. V\'azquez, \emph{Fractional Nonlinear Degenerate Diffusion Equations on Bounded Domains Part II. Positivity, Boundary behaviour and Harnack inequalities. }In preparation (2016).

\bibitem{BV-ADV}  {\rm M.  Bonforte, J.~L. V\'azquez. } {\it Positivity, local smoothing, and Harnack inequalities for very fast diffusion equations}, \textrm{Adv. Math.} \bf 223 \rm (2010), 529--578.

\bibitem{BV2012} M. Bonforte, J.~L. V\'azquez. \textit{Quantitative Local and Global A Priori Estimates for Fractional Nonlinear Diffusion Equations}, Adv. Math. \textbf{250} (2014), 242--284. 

\bibitem{CaFr79b}L. A. Caffarelli, A. Friedman. \textit{Continuity of the density of a gas flow in a porous medium. }Trans. Amer. Math. Soc. 252 (1979), 99–113.

\bibitem{CaFr80}L. A. Caffarelli, A. Friedman. \textit{Regularity of the free boundary of a gas flow in an n-dimensional porous medium. }Indiana Univ. Math. J. \textbf{29} (1980), 361--391.

\bibitem{CS2011}L. A. Caffarelli, L. Silvestre. \textit{Regularity results for nonlocal equations by approximation. }
Arch. Rat. Mech. Anal. \textbf{200} (2011), no. 1, 59--88

\bibitem{CV2011} {\rm L. A. Caffarelli, J.L. Vazquez,}, \emph{Nonlinear porous medium flow with fractional potential pressure}.
{Arch. Rat. Mech. Anal.} \textbf{202} (2011) 537–-565.

\bibitem{CV2010}  {\rm L. A. Caffarelli, J.L. Vazquez,}, \emph{Asymptotic behaviour of a porous medium equation with fractional
diffusion}. {Disc. Cont. Dyn. Syst.} \textbf{29} (2011), no.~4, 1393--1404.

\bibitem{CaVaWo87}L. A. Caffarelli, J.L. Vazquez, N. Wolansky. \textit{Lipschitz continuity of solutions and interfaces of the n-dimensional porous medium equation}, Ind. Univ. Math. J. \textbf{36} (1987), 373--401.
%
\bibitem{CaWo90}L. A. Caffarelli, N. Wolansky. \textit{$C^{1,\alpha}$-regularity of the free boundary for the n-dimensional porous media equation}, Comm. Pure and Appl. Math. \textbf{43} (1990), 885--902.

\bibitem{counterexample} H. Chang-Lara, G. Davila, \emph{Regularity for solutions of nonlocal parabolic equations II}, J. Differential Equations 256 (2014), 130-156.

\bibitem{CKriv} H. Chang-Lara, D. Kriventsov, \emph{Further time regularity for non-local, fully non-linear parabolic equations}, Comm. Pure Appl. Math., to appear.

\bibitem{Chen-Song} Z.-Q. Chen, R. Song, \emph{Estimates on Green functions and Poisson kernels for symmetric stable processes}, Math. Ann. 312 (1998), 465-501.

\bibitem{finance2} R. Cont, P. Tankov, \emph{Financial Modelling With Jump Processes}, Chapman \& Hall, Boca Raton, FL, 2004.

\bibitem{CP-JFA} M. Crandall, M. Pierre, \textit{Regularizing Effects for $u_t=A\varphi(u)$ in $\LL^1$}, J. Funct. Anal. \textbf{45}, (1982), 194-212


\bibitem{DaHaLe01}P. Daskalopoulos, R. Hamilton, K. Lee, \textit{All time $C^\infty$-regularity of the interface in degenerate diffusion: a geometric approach. }Duke Math. J. \textbf{108} (2001), no. 2, 295--327.
%
\bibitem{DaHa98}P. Daskalopoulos, R. Hamilton, \textit{Regularity of the free boundary for the porous medium equation. }J. Amer. Math. Soc. \textbf{11} (1998), no. 4, 899--965.

\bibitem{DaskaBook}P. Daskalopoulos, C. E. Kenig, \textsl{``Degenerate diffusions. Initial value problems and local regularity theory''. }EMS Tracts in Mathematics, \textbf{1}. European Mathematical Society (EMS), Zürich, 2007. x+198 pp. ISBN: 978-3-03719-033-3

\bibitem{DiB86}E. DiBenedetto. \textit{On the local behaviour of solutions of degenerate parabolic equations with measurable coefficients. }Ann. Scuola Norm. Sup. Pisa Cl. Sci. (4) \textbf{13} (1986), no. 3, 487--535.
%
\bibitem{DiB88}E. DiBenedetto. \textit{Intrinsic Harnack type inequalities for solutions of certain degenerate parabolic equations. }Arch. Rational Mech. Anal. \textbf{100} (1988), no. 2, 129--147.

\bibitem{DiBook}E. DiBenedetto. \textsl{``Degenerate parabolic equations''}. Universitext. Springer-Verlag, New York, 1993. xvi+387 pp. ISBN: 0-387-94020-0
%
\bibitem{DGVacta} E. DiBenedetto, U. Gianazza, V. Vespri. \textit{Harnack estimates for quasi-linear degenerate parabolic differential equations. } Acta Math. \textbf{200} (2008), no. 2, 181--209.
%
\bibitem{DGVbook} E. DiBenedetto, U. Gianazza, V. Vespri. \textsl{``Harnack'\,s inequality for degenerate and singular parabolic equations''}\,,  Springer Monographs in Mathematics, Springer 2011.

\bibitem{DL1972} {\rm G. Duvaut, J.-L. Lions, }{\em ``Les In{\'e}quations en Mechanique et en Physique''}, Travaux et Recherches Mathématiques, No. 21. Dunod, Paris, 1972. xx+387 pp.

\bibitem{FGS86} E. B. Fabes, N. Garofalo, S. Salsa. \textit{``A backward Harnack inequality and Fatou theorem for nonnegative solutions of parabolic equations''}, Illinois J. Math. \textbf{30} n. 4 (1986), 536--565.

\bibitem{FK} M. Felsinger, M. Kassmann. \textit{Local regularity for parabolic nonlocal operators}, Comm. Partial Differential Equations \textbf{38} (2013), 1539-1573.

\bibitem{GiPe81}B. H. Gilding, L. A. Peletier, \textit{Continuity of solutions of the porous media equation. }Ann. Scuola Norm. Sup. Pisa Cl. Sci. (4) \textbf{8} (1981), no. 4, 659--675.

\bibitem{Grub1} G. Grubb, \textit{Fractional Laplacians on domains, a development of H\"ormander's theory of mu-transmission pseudodifferential operators}, Adv. in Math. \textbf{268} (2015), 478--528.

\bibitem{nature} N. E. Humphries et al., \emph{Environmental context explains Levy and Brownian movement patterns of marine predators}, Nature 465 (2010), 1066-1069.

\bibitem{Jara0} {\rm M. Jara,} \textit{Nonequilibrium scaling limit for a tagged particle in the simple exclusion process with long jumps, } \textrm{Comm. Pure Appl. Math.}, {\bf 62} (2009),  198--214.

\bibitem{Jara1} {\rm M. Jara}, \textit{Hydrodynamic Limit Of Particle Systems With Long Jumps}\,, {\rm Preprint, }(2009). \\ \texttt{http://arxiv.org/abs/0805.1326v2}

\bibitem{JKOlla} {\rm M. D. Jara, T. Komorowski, S.  Olla,} \textit{Limit theorems for additive functionals of a Markov chain}. \textrm{Ann. Appl. Probab.} \textbf{19} (2009), no. 6, 2270–-2300.

\bibitem{Jara2} {\rm M. Jara, C. Landim,  S. Sethuraman, } \textit{Nonequilibrium fluctuations for a tagged particle in mean-zero one-dimensional zero-range processes, }\textrm{ Probab. Theory Relat. Fields }{\bf 145 }(2009), 565--590.

\bibitem{JX} T. Jin, J. Xiong \textit{Schauder estimates for solutions of linear parabolic integro-differential equations}, Disc. Contin. Dyn. Syst. 35 (2015), 5977-5998.

\bibitem{Kim-Coeff} K.-Y. Kim, P. Kim, \emph{Two-sided estimates for the transition densities of symmetric Markov processes dominated by stable-like processes in $C^{1,\eta}$ open sets}. Stochastic Processes and their Applications 124 (2014) 3055--3083

\bibitem{Kul} T. Kulczycki. \textit{Properties of Green function of symmetric stable processes. }Probab. Math. Statist. \textbf{17} (1997), 339--364.

\bibitem{phys} N. Laskin, \emph{Fractional quantum mechanics and L\'evy path integrals}, Physics Letters 268 (2000), 298-304.

\bibitem{LMT2003} E. K. Lenzi, R. S. Mendes, C. Tsallis, \emph{Crossover in diffusion equation: Anomalous and normal behaviors}, {Physical Review E }\textbf{67}, 031104 (2003).

\bibitem{Levy} P. L\'evy, \emph{Th\'eorie de l'addition des variables al\'eatoires}, Gauthier-Villars, Paris, 1937.

\bibitem{phys2} R. Metzler, J. Klafter, \emph{The random walk's guide to anomalous diffusion: a fractional dynamics approach}, Phys. Rep. 339 (2000), 1-77.

\bibitem{DPQR} {\rm A.~de Pablo, F. Quir\'os, A. Rodr\'iguez. }\textit{Nonlocal filtration equations with rough kernels.} Preprint (2015). Arxiv: \texttt{http://arxiv.org/abs/1509.09143}

\bibitem{DPQRV1} {\rm A.~de Pablo, F. Quir\'os, A. Rodr\'iguez, J. L. V\'azquez. }\textit{A fractional porous medium equation} Adv. Math. \textbf{226} (2011), no. 2, 1378--1409.

\bibitem{DPQRV2}{\rm A.~de Pablo, F. Quir\'os, A. Rodr\'iguez, J. L. V\'azquez. } \textit{A general fractional porous medium equation},  Comm. Pure Appl. Math. {\bf 65} (2012), 1242--1284.

\bibitem{nature2} A. M. Reynolds, C. J. Rhodes, \emph{The L\'evy flight paradigm: Random search patterns and mechanisms}, Ecology 90 (2009), 877-887.

\bibitem{R-survey} X. Ros-Oton. \emph{Nonlocal elliptic equations in bounded domains: a survey}, {Publ. Mat.} 60 (2016), 3-26.

\bibitem{RS-Dir} X. Ros-Oton, J. Serra. \emph{The Dirichlet problem for the fractional Laplacian: regularity up to the boundary}, {J.~Math. Pures Appl.} 101 (2014), 275-302.

\bibitem{finance} W. Schoutens, \emph{L\'evy Processes in Finance: Pricing Financial Derivatives}, Wiley, 2003.

\bibitem{SdTV1} D. Stan, F. del Teso, J. L. Vazquez, \textit{Transformations of self-similar solutions for porous medium equations of fractional type. }Nonlinear Anal. \textbf{119} (2015), 62--73.

\bibitem{SdTV2} D. Stan, F. del Teso, J. L. Vazquez, \textit{Finite and infinite speed of propagation for porous medium equations with fractional pressure. }C. R. Math. Acad. Sci. Paris \textbf{352} (2014), no. 2, 123--128.

\bibitem{SdTV3} D. Stan, F. del Teso, J. L. Vazquez, \textit{Finite and infinite speed of propagation for porous medium equations with nonlocal pressure, }To appear in J. Differential Equations (2015).

\bibitem{SV2013} {  D. Stan, J.~L. V{\'a}zquez, } \emph{The Fisher-KPP equation with nonlinear fractional diffusion}, SIAM J. Math. Anal. \textbf{46} (5), 3241--3276

\bibitem{VDPQR} {\rm  J.~L.~V\'azquez,  A. de Pablo,  F. Quir\'os, A. Rodr\'iguez. }{\em Classical solutions and higher regularity for nonlinear fractional diffusion equations}, To appear in  J. Eur. Math. Soc. (2015).

\bibitem{VazLN} {\rm  J.~L. V{\'a}zquez. }{\it Smoothing and decay estimates for nonlinear diffusion equations}, vol.~33 of Oxford Lecture Notes in Maths. and its Applications, Oxford Univ. Press, 2006.

\bibitem{VazBook}{\rm  J.~L. V{\'a}zquez. } \emph{The Porous Medium Equation. Mathematical Theory}, vol.~Oxford Mathematical Monographs, Oxford University Press, Oxford, 2007.

\bibitem{VazAbel}{\rm  J.~L. V{\'a}zquez. }{\sl Nonlinear Diffusion with Fractional Laplacian Operators.}  in ``Nonlinear partial differential equations: the Abel Symposium 2010'', Holden, Helge  \& Karlsen, Kenneth H. eds., Springer, 2012. Pp. 271--298.

\bibitem{VazSurvey}{\rm  J.~L. V{\'a}zquez. }  {\sl Recent progress in the theory  of Nonlinear Diffusion with  Fractional Laplacian Operators}. In ``Nonlinear elliptic and parabolic differential equations'', Disc. Cont. Dyn. Syst. - S \ {\bf 7,} no. 4 (2014), 857--885.

\end{thebibliography}
\end{document}